\documentclass[twoside,11pt,reqno]{amsart}
\usepackage{amsmath,amssymb,amscd,mathrsfs,epic,wasysym,latexsym,tikz,mathrsfs,cite,hyperref}
\usepackage{pb-diagram}
\usepackage[matrix,arrow]{xy}

\makeatletter

\numberwithin{equation}{section}

\newtheorem{Proposition}[equation]{Proposition}
\newtheorem{Lemma}[equation]{Lemma}
\newtheorem{Theorem}[equation]{Theorem}
\newtheorem{Corollary}[equation]{Corollary}
\theoremstyle{definition}  
\newtheorem{Definition}[equation]{Definition}
\newtheorem{Remark}[equation]{Remark}

\newtheorem{Conjecture}[equation]{Conjecture}

\let\<\langle
\let\>\rangle

\newcommand\Comment[2][\relax]{\space\par\medskip\noindent%
   \fbox{\begin{minipage}{\textwidth}\textbf{Comment\ifx\relax#1\else---#1\fi}\newline%
        #2\end{minipage}}\medskip
}

\def\bi{\text{\boldmath$i$}}
\def\bj{\text{\boldmath$j$}}

\def\b1{\text{\boldmath$1$}}

\def\pmod#1{\text{ }(\text{\rm mod } #1)\,}

\newcommand{\Hom}{\operatorname{Hom}}
\newcommand{\Ext}{\operatorname{Ext}}
\newcommand{\EXT}{\operatorname{Ext}}

\newcommand{\End}{\operatorname{End}}

\newcommand{\im}{\operatorname{im}}

\newcommand{\soc}{\operatorname{soc}\,}
\newcommand{\head}{\operatorname{head}}

\newcommand{\noncusp}{\operatorname{nsc}}

\newcommand{\Z}{\mathbb{Z}}

\def\eps{{\varepsilon}}
\def\phi{{\varphi}}

\newcommand{\catC}{{\mathbf C}}

\newcommand{\funQ}{{\mathcal Q}}

\newcommand{\funR}{{\mathcal R}}
\newcommand{\funE}{{\mathcal E}}
\newcommand{\funO}{{\mathcal O}}

\newcommand{\ga}{\gamma}
\newcommand{\Ga}{\Gamma}
\newcommand{\la}{\lambda}
\newcommand{\La}{\Lambda}
\newcommand{\al}{\alpha}
\newcommand{\be}{\beta}

\def\Si{\mathfrak{S}}
\newcommand{\si}{\sigma}

\newcommand{\de}{\delta}
\newcommand{\De}{\Delta}

\newcommand{\Ind}{{\mathrm {Ind}}}
\newcommand{\Coind}{{\mathrm {Coind}}}

\newcommand{\rad}{{\mathrm {rad}\,}}

\newcommand{\Res}{{\mathrm {Res}}}
\newcommand{\Ann}{{\mathrm {Ann}}}

\newcommand{\Q}{{\mathbb Q}}

\newcommand{\D}{{\mathscr D}}

\renewcommand{\mod}{\bmod \,}

\def\h{{\mathfrak h}}
\def\g{{\mathfrak g}}

\def\Par{{\mathscr P}}
\def\ula{{\underline{\lambda}}}
\def\umu{{\underline{\mu}}}

\def\b{\mathfrak{b}}
\def\k{\Bbbk}

\def\height{{\operatorname{ht}}}

\def\op{{\mathrm{op}}}
\def\re{{\mathrm{re}}}
\def\im{{\mathrm{im}\,}}

\def\mod#1{#1\!\operatorname{-mod}}

\renewcommand\O{\mathcal O}
\def\HCI{I}
\def\HCR{{}^*I}
\def\iso{\stackrel{\sim}{\longrightarrow}}

\def\Mde{M}

\def\HOM{\operatorname{Hom}}

\def\CH{{\operatorname{ch}_q\,}}
\def\DIM{{\operatorname{dim}_q\,}}

\def\words{I}

\def\Car{{\tt C}}
\def\cc{{\tt c}}

\def\ImS{{\mathscr S}}

{\catcode`\|=\active
  \gdef\set#1{\mathinner{\lbrace\,{\mathcode`\|"8000%
  \let|\midvert #1}\,\rbrace}}
}
\def\midvert{\egroup\mid\bgroup}

\colorlet{darkgreen}{green!50!black}
\tikzset{dots/.style={very thick,loosely dotted},
         greendot/.style={fill,circle,color=darkgreen,inner sep=1.5pt,outer sep=0},
         blackdot/.style={fill,circle,color=black,inner sep=1.5pt,outer sep=0},
         graydot/.style={fill,circle,color=gray,inner sep=1.1pt,outer sep=0}
}
\def\greendot(#1,#2){\node[greendot] at(#1,#2){}}
\def\blackdot(#1,#2){\node[blackdot] at(#1,#2){}}
\def\graydot(#1,#2){\node[graydot] at(#1,#2){}}

\newenvironment{braid}{
  \begin{tikzpicture}[baseline=6mm,black,line width=1pt, scale=0.32,
                      draw/.append style={rounded corners},
                      every node/.append style={font=\fontsize{5}{5}\selectfont}]%
  }{\end{tikzpicture}
}

\def\Grid(#1,#2){
  \draw[very thin,gray,step=2mm] (0,0)grid(#1,#2);
  \draw[very thin,darkgreen,step=10mm] (0,0)grid(#1,#2);
}

\newcommand\Tableau[2][\relax]{
  \begin{tikzpicture}[scale=0.5,draw/.append style={thick,black}]
    \ifx\relax#1\relax%
    \else 
      \foreach\box in {#1} { \filldraw[blue!30]\box+(-.5,-.5)rectangle++(.5,.5); }
    \fi
    \newcount\row\newcount\col
    \row=0
    \foreach \Row in {#2} {
       \col=1
       \foreach\k in \Row {
          \draw(\the\col,\the\row)+(-.5,-.5)rectangle++(.5,.5);
          \draw(\the\col,\the\row)node{\k};
          \global\advance\col by 1
       }
       \global\advance\row by -1
    }
  \end{tikzpicture}
}

\newcommand\YoungDiagram[2][\relax]{
  \begin{tikzpicture}[scale=0.5,draw/.append style={thick,black}]
    \ifx\relax#1\relax%
    \else 
    \foreach\box in {#1} {
      \filldraw[blue!30]\box rectangle ++(1,1);
    }
    \fi
    \newcount\row
    \row=0
    \foreach \col in {#2} {
       \draw(1,\the\row)grid ++(\col,1);
       \global\advance\row by -1
    }
  \end{tikzpicture}
}

\begin{document}

\title[Stratifying KLR algebras of affine ADE types]{{\bf Stratifying KLR algebras of affine ADE types}}

\author{\sc Alexander Kleshchev}
\address{Department of Mathematics\\ University of Oregon\\
Eugene\\ OR 97403, USA}
\email{klesh@uoregon.edu}

\author{\sc Robert Muth}
\address{Department of Mathematics\\ Tarleton State University\\
Stephenville \\ TX 76402, USA}
\email{robmuth@gmail.com}

\subjclass[2000]{20C08, 20C30, 05E10}

\thanks{
Supported by the NSF grant DMS-1161094, Max-Planck-Institut and Fulbright Foundation.}

\begin{abstract}
We generalize imaginary Howe duality for KLR algebras of affine $ADE$ types, developed in our previous paper, from balanced to arbitrary convex preorders. Under the assumption that the characteristic of the ground field is greater than some explicit bound, we prove that these KLR algebras are properly stratified. 
\end{abstract}

\maketitle

\centerline
{\em Dedicated to the memory of Professor J.A. Green.}

\vspace{5mm}

\section{Introduction}
There are different ways to study representation theory of Khovanov-Lauda-Rouquier (KLR) algebras of \cite{KL1} and \cite{Ro}. For example, crystal-theoretic methods are used in \cite{LV} to classify simple modules over these algebras, while connections to dual canonical bases are explored in \cite{VV}, \cite{Ro2}, \cite{KKP}. In finite and affine types $A$, there are combinatorial methods coming from the theory of cyclotomic Hecke algebras \cite{BKcyc}, \cite{BKW}, \cite{BKgrdec}, \cite{HM}.

An approach through some versions of standard modules, indicating a presence of a stratification structure on the algebras, has been pursued in \cite{KRbz}, \cite{BKOP}, \cite{Kato}, \cite{McN}, \cite{BKM}, \cite{Kcusp}, \cite{TW}, \cite{KM}, \cite{McNAff}. 
In this paper we work with the KLR algebras $R_\theta$ of untwisted affine $ADE$ types over an arbitrary field $k$ of characteristic $p\geq 0$. The goal of this paper is two-fold. 

First, we generalize much of the theory developed in \cite{KM} from balanced to arbitrary convex preorders. In particular, we obtain an analogue of the imaginary Howe duality theory in complete generality (for untwisted affine $ADE$ types). In addition to the methods used in \cite{KM}, we need two new ingredients: work of McNamara \cite{McNAff}, which gives the desired result in characteristic zero, and reduction modulo $p$. 

Under the assumption  $p=0$, it is proved in \cite{McNAff} that $R_\theta$ is properly stratified. Informally, this means that the category $\mod{R_\theta}$ of finitely generated graded $R_\theta$-modules is stratified  by the categories $\mod{B_\xi}$ for much simpler algebras $B_\xi$. Our second goal is to apply reduction modulo $p$ arguments to generalize this result to the case where $p$ is greater than some explicit bound. 

Description of the algebras $B_\xi$ appearing in the stratification of  $\mod{R_\theta}$ is easily reduced to the semicuspidal cases, which split into  real and imaginary subcases. In the real case we have $B_{n\al}\cong k[z_1,\dots,z_n]^{\Si_n}$,  
but in the imaginary case $B_{n\de}$ is not so easy to understand. 
In the sequel paper \cite{KMZZ}, we reveal a connection between $B_{n\de}$ and affine zigzag algebras related to the {\em `finite zigzag algebras'} of \cite{HK} of the underlying finite Lie type.

We now describe the contents of the paper in more detail. In Section~\ref{SStrat}, we follow \cite{Kdonkin} to recall the formalism of properly stratified graded algebras, which goes back to \cite{CPS}. A choice of a partial preorder on the set $\Pi$ of labels of simple modules allows one to define the standard modules $\{\De(\pi)\mid \pi\in\Pi\}$. The requirement is then that the projective covers $P(\pi)$ have  standard filtrations of special kind and that standard modules are finitely generated flat over their  endomorphism algebras, see Definition~\ref{DStCat}. 

In Section~\ref{SLie}, we fix Lie theoretic notation, in particular, the simple roots $\{\al_i\mid i\in I\}$, where $I=\{0,1,\dots,l\}$ with $0$ corresponding to the affine simple root, and the positive part of the root lattice  $Q_+=\sum_{i\in I}\Z_{\geq 0}\cdot \al_i$. We define a convex preorder $\preceq$ on the set $\Phi_+$ of positive roots  and discuss its properties. We define Kostant partitions and root partitions. Kostant partitions of $\theta$ are going to label the strata in $\mod{R_\theta}$, while a more refined notion of a root partition of $\theta$ will be needed to label the irreducible modules in $\mod{R_\theta}$. 

In Section~\ref{SKLR}, we define KLR algebras $R_\theta$ for $\theta\in Q_+$ and discuss their basic representation theory. Of particular importance are induction and restriction functors $\Ind_{\underline{\theta}}$ and $\Res_{\underline{\theta}}$. We also use the notation `$\circ$' when inducing from outer tensor products. Finally, we discuss reduction modulo $p$. For this we use the $p$-modular system $(F,\O,K)$ with $F={\mathbb F}_p$, $\O=\Z_p$ and $K=\Q_p$. 

In Section~\ref{SSCSM}, we recall the notion of semicuspdial representations of $R_{n\al}$, where $\al$ is an indivisible positive root. The semicuspidal algebra $C_{n\al}$ is defined so that the category of finitely generated semicuspidal $R_{n\al}$-modules is equivalent to $\mod{C_{n\al}}$. The main goal of the section is to develop a general theory of imaginary Howe duality, which allows us to classify irreducible $C_{n\de}$-modules. The main results  are similar to those contained in \cite{KM} and \cite{McNAff}, except that in \cite{KM} the convex order is assumed to be balanced, while in \cite{McNAff} the ground field is assumed to be of characteristic zero. 

In \S\ref{SSMinusc}, we construct the minuscule imaginary modules $L_{\de,i}$, one for each $i\in I':=\{1,\dots,l\}$, and prove that reduction modulo $p$ of $L_{\de,i,K}$ is $L_{\de,i,F}$. Next, we fix $i\in I'$ and review the action of the symmetric group $\Si_n$ on the imaginary tensor space $M_{n,i}:= L_{\de,i}^{\circ n}$. We normalize this action using  Gelfand-Graev modules $\Ga_{n,i}$ defined in \S\ref{SSISW}. 
The action of the symmetric group allows us to define the imaginary divided powers $Z_{n,i}:=M_{n,i}^{\Si_n}$. 

Fix $h\geq n$ and denote by $\La(h,n)$ the set of compositions of $n$ with $h$ parts. For  $\nu:=(n_1,\dots,n_h)\in\La(h,n)$, we write $n=|\nu|$ and define $Z^\nu_i:=Z_{n_1,i}\circ\dots\circ Z_{n_h,i}$. We set 
\begin{equation}\label{EBigZ}
Z(h,n):=\bigoplus_{\nu(1),\dots,\nu(l)} Z^{\nu(1)}_1\circ \dots\circ Z^{\nu(l)}_l,
\end{equation}
where the sum is over all compositions $\nu(1),\dots,\nu(l)$ with $h$ parts such that $|\nu(1)|+\dots+|\nu(l)|=n$. Finally, we define the imaginary Schur algebra $\ImS(h,n)$ as $R_{n\de}/\Ann_{R_{n\de}}(Z(h,n))$. 

\vspace{2 mm}
\noindent
{\bf Theorem 1.}
{\em
For $h\geq n$, the $\ImS(h,n)$-module $Z(h,n)$ is a projective generator in $\mod{\ImS(h,n)}$, and 
$$\End_{\ImS(h,n)}(Z(h,n))^\op\cong \bigoplus_{(n_1,\dots,n_l)}S(h,n_1)\otimes \dots\otimes S(h,n_l),
$$ where the sum is over all compositions $(n_1,\dots,n_l)$ of $n$ with $l$ parts, and $S(h,n)$ denotes the classical Schur algebra. 
}
\vspace{2 mm}

Denote 
\begin{equation}\label{EBigSchur}
{\bf S}(h,n):= \bigoplus_{(n_1,\dots,n_l)}S(h,n_1)\otimes \dots\otimes S(h,n_l), 
\end{equation}
where the sum is as in the theorem. By the classical theory, the irreducible ${\bf S}(h,n)$-modules are labeled by the set  $\Par_n$ of $l$-multipartitions of $n$. 
By Theorem 1, $\ImS(h,n)$ is Morita equivalent to ${\bf S}(h,n)$,
so we get the irreducible $\ImS(h,n)$-modules $\{L(\umu)\mid\umu\in\Par_n\}$. Lifting to $R_{n\de}$, these can be considered as $R_{n\de}$-modules. 

\vspace{2 mm}
\noindent
{\bf Theorem 2.}
{\em
The set $\{L(\umu)\mid\umu\in\Par_n\}$ is a complete and irredundant set of semicuspidal irreducible $R_{n\de}$-modules up to isomorphism and degree shift. Moreover, if $p>n$, then reduction modulo $p$ of $L(\umu)_K$ is $L(\umu)_F$. 
}
\vspace{2 mm}

Theorem 1 and the first statement of Theorem 2 are proved in \S\ref{SST1}. 

Section~\ref{SStratKLR} is on the stratification of $R_\theta$. 
Projective indecomposable modules in $\mod{C_{n\al}}$ are used to define standard modules for $R_\theta$. We show that our  definitions, which use parabolic induction of semicuspidal representations, agree with a general categorical definition. In Theorem~\ref{TFlatness}, we verify the flatness condition in the definition of properly stratified algebras. To verify the standard filtration condition we need a certain $\Ext$ result, which is proved in  Theorem~\ref{TExt} following McNamara's argument in  \cite{McNAff}.  

\vspace{2 mm}
\noindent
{\bf Theorem 3.}
{\em
Let $\theta=\sum_{i\in I}n_i\al_i\in Q_+$ and assume that $p>\min\{n_i\mid i\in I\}$. For any convex preorder on $\Phi_+$, the algebra $R_\al$ is properly stratified. 
Moreover, for any $\pi\in\Pi(\theta)$, the standard module $\De(\pi)_F$ is a reduction modulo $p$ of $\De(\pi)_K$; in particular, $\De(\pi)_F$ and $\De(\pi)_K$ have the same formal characters. The same is true for proper standard modules $\bar\De(\pi)_F$ and $\bar\De(\pi)_K$. 
}
\vspace{2 mm}

The results on reduction modulo $p$ cited in Theorems 2 and 3 are proved in Section~\ref{SRed}. 

\subsection*{Acknowledgement} We are grateful to Peter McNamara for useful explanations of \cite{McNAff}. 

\section{Stratification}\label{SStrat}
Throughout the paper, unless otherwise stated, $k$ is an arbitrary  field of characteristic $p\geq 0$.  Let $q$ be a variable, and $\Z((q))$ be the ring of Laurent series. The quantum integers $[n]=(q^n-q^{-n})/(q-q^{-1})$ as well as  expressions like $[n]!:=[1][2]\dots[n]$ and $1/(1-q^2)$ are always interpreted as elements of $\Z((q))$.

In this section, we mainly follow \cite{Kdonkin}. All notions we consider, such as algebras, modules, ideals, etc., are assumed to be ($\Z$-)graded. 

\subsection{Graded algebras} 
If $H$ is a Noetherian (graded) $k$-algebra, we denote by $\mod{H}$ the category of finitely generated graded left $H$-modules. The morphisms in this category are all homogeneous degree zero  $H$-homomorphisms, which we denote $\hom_{H}(-,-)$. 
For  $V\in\mod{H}$, let $q^d V$ denote its grading shift by $d$,  so if $V_m$ is the degree $m$ component of $V$, then $(q^dV)_m= V_{m-d}.$ For a polynomial $a=\sum_{d}a_dq^d\in\Z[q,q^{-1}]$ with non-negative coefficients, we set $a V:=\bigoplus_d(q^d V)^{\oplus a_d}$.

For $U,V\in \mod{H}$,
we set 
$$\HOM_H(U, V):=\bigoplus_{d \in \Z} \HOM_H(U, V)_d,$$ 
where
$
\HOM_H(U, V)_d := \hom_H(q^d U, V) .
$
We define $\operatorname{ext}^m_H(U,V)$ and
$\EXT^m_H(U,V)$ similarly. 
Since $U$ is finitely generated, $\HOM_H(U,V)$ can be identified with the set of all $H$-module homomorphisms  ignoring the gradings. A similar result holds for $\EXT^m_H(U,V)$, since $U$ has a resolution by finitely generated projective modules. Given $V,W\in\mod{H}$, we write $V\simeq W$ to indicate that $V\cong q^n W$ for some $n\in\Z$.

A vector space $V$ is called {\em Laurentian} if its graded components $V_m$ are finite dimensional and $V_m=0$ for $m\ll0$. Then the graded dimension $\DIM V$ is a Laurent series.  
An algebra is called {\em Laurentian} if it is so as a vector space.  
In this case all irreducible $H$-modules are finite dimensional, there are only finitely many irreducible $H$-modules up to isomorphism and degree shift, and 
every finitely generated $H$-module has a projective cover, see  \cite[Lemma 2.2]{Ksing}. 

Let $H$ be a Laurentian algebra. We fix a complete irredundant set $\{L(\pi)\mid \pi\in\Pi\}$ of irreducible $H$-modules up to isomorphism and degree shift. 
By above, the set $\Pi$ is finite. 
For each $\pi\in\Pi$, we also fix a projective cover $P(\pi)$ of $L(\pi)$. 
Let 
\begin{equation}\label{EM}
M(\pi):=\rad P(\pi) \qquad(\pi\in \Pi), 
\end{equation}
so that $P(\pi)/M(\pi)\cong L(\pi)$.  

From now on we assume in addition that $H$ is {\em Schurian}, i.e. $\End(L(\pi))\cong k$ for all $\pi$. 
For any $V\in\mod{H}$ and $\pi\in\Pi$, the  {\em composition multiplicity} of $L(\pi)$ in $V$ is defined as $[V:L(\pi)]_q:=\DIM\Hom(P(\pi),V)\in\Z((q))$. 

We record for future references:

\begin{Lemma} \label{LSchubert} {\rm \cite[Lemma 4.3.1]{BDK}} 
Let $A$ be a (graded) algebra and $0 \to Z \to P \to M \to 0$ be a short exact sequence of (graded) $A$-modules with $P$ (graded) projective. If every (degree zero) $A$-module homomorphism from $P$ to $M$ annihilates $Z$, then $M$ is a (graded) projective  $A/\Ann_A(M)$-module.
\end{Lemma}

\subsection{Standard objects and stratification}\label{SSStObStr}
We continue with the notation of the previous subsection. 
Let $\Sigma$ be a subset of $\Pi$. An object $X$ of the category $\catC:=\mod{H}$ {\em belongs to $\Sigma$} if $[X:L(\si)]_q\neq 0$ implies $\si\in\Sigma$. 
Let $\catC(\Sigma)$ be the full subcategory of $\catC$ consisting of all objects which belong to $\Sigma$. 
The natural inclusion $\iota_\Sigma:\catC(\Sigma)\to\catC$ 
has left adjoint $\funQ^\Sigma:\catC\to\catC(\Sigma)$ with  $\funQ^\Sigma(V)=V/\funO^\Sigma(V)$, where $\funO^\Sigma(V)$ is the unique minimal submodule among submodules $U\subseteq V$ such that $V/U$ belongs to $\Sigma$. Let also $\funO_\Sigma(V)$ be the unique maximal subobject of $V$ which belongs to $\Sigma$. 

Applying $\funO^\Sigma$ to the left regular module $H$ yields a (two-sided) ideal  $\funO^\Sigma(H)\unlhd H$. 
By \cite[Lemma 3.12]{Kdonkin}, for  $V\in\mod{H}$, we have $\funO^\Sigma(H)V=\funO^\Sigma(V)$. 
Set $
H(\Sigma):=H/\funO^\Sigma(H).
$
Then we can regard $\funQ^\Sigma(V)$ as an $H(\Sigma)$-module and identify $\catC(\Sigma)$ and $\mod{H(\Sigma)}$. In this way,   $\funQ^\Sigma$ becomes a functor
$
\funQ^\Sigma:\mod{H}\to\mod{H(\Sigma)}.
$

We now suppose that there is a fixed surjection 
\begin{equation}\label{ERho}
\varrho: \Pi\to \Xi
\end{equation}
for some set $\Xi$ endowed with a    
{\em partial order}\, $\leq$. We then have a partial preorder $\leq$ on $\Pi$ with $\pi\leq\si$ if and only if $\varrho(\pi)\leq\varrho(\si)$. 
For $\pi\in \Pi$ and $\xi\in \Xi$ we set 
\begin{align*}
&\Pi_{< \pi}:=\{\si\in\Pi\mid \si< \pi\},\ \Pi_{\leq \pi}:=\{\si\in\Pi\mid \si\leq \pi\},
\\
&\Pi_{< \xi}:=\{\si\in\Pi\mid \varrho(\si)< \xi\},\ \Pi_{\leq \xi}:=\{\si\in\Pi\mid \varrho(\si)\leq \xi\},
\end{align*}
and write 
$
\funO^{\leq \pi}:=\funO^{\Pi_{\leq \pi}}$, $\funO_{<\xi}:=\funO_{\Pi_{< \xi}}$, $ 
\funQ^{< \pi}:=\funQ^{\Pi_{< \pi}} 
$, 
$H_{\leq\xi}:=H(\Pi_{\leq \xi})$, etc.

Recalling (\ref{EM}), we define for all $\pi\in\Pi$: 
$$
K(\pi):=\funO^{\leq \pi}(P(\pi))=\funO^{\leq \pi}(M(\pi)),\qquad \bar K(\pi):=\funO^{< \pi}(M(\pi)), 
$$
and 
\begin{equation}\label{EStandard}
\De(\pi):=\funQ^{\leq \pi}(P(\pi))=P(\pi)/K(\pi),\qquad \bar\De(\pi):=P(\pi)/\bar K(\pi).
\end{equation}
Note that $\bar K(\pi)\supseteq K(\pi)$, so $\bar \De(\pi)$ is naturally a quotient of $\De(\pi)$. Moreover, $\head \De(\pi)\cong \head \bar\De(\pi)\cong L(\pi)$. 
We call the modules $\De(\pi)$ {\em standard} and the modules $\bar\De(\pi)$ {\em proper standard}. 
By \cite[Lemma~3.10]{Kdonkin}, $\De(\pi)$ is the projective cover of $L(\pi)$ in the category $\catC_{\leq \pi}$. 
For  $V\in\catC$, a finite {\em $\De$-filtration} (or a {\em standard filtration}) is a filtration
$
V=V_0\supseteq V_1\supseteq \dots \supseteq V_N=0
$ such that each for $0\leq n<N$ we have $V_n/V_{n+1}\simeq\De(\pi)$ for some $\pi\in \Pi$.

Let $\xi\in\Xi$. Then $\catC_{<\xi}$ is a  
Serre subcategory of $\catC_{\leq \xi}$, and the quotient category 
\begin{equation}\label{ECXi}
\catC_\xi:=\catC_{\leq \xi}/\catC_{< \xi},
\end{equation} 
is called the {\em $\xi$-stratum}. 
Up to isomorphism and degree shift, $\{L(\pi)\mid\varrho(\pi)=\xi\}$ is a complete family of simple objects in $\catC_\xi$, and 
$
P_\xi(\pi):=\De(\pi)/\funO_{<\xi}(\De(\pi))
$ 
is the projective cover of $L(\pi)$ in $\catC_\xi$.  
Finally, setting  
\begin{equation}\label{EDeXi}
\De(\xi):=\bigoplus_{\pi\in\varrho^{-1}(\xi)}\De(\pi)\quad \text{and}\quad B_\xi:=\End_H(\De(\xi))^{\op},
\end{equation}
by \cite[Corollary 4.9]{Kdonkin}, the stratum category $\catC_\xi$ is graded equivalent to $\mod{B_\xi}$.

We have a natural exact projection functor $\funR_\xi:\catC_{\leq \xi}\to\catC_\xi$. If we identify $\catC_\xi$ and $\mod{B_\xi}$, the functor $R_\xi$ becomes 
$$
\funR_\xi=\Hom_{H_{\leq\xi}}(\De(\xi),-):\mod{H_{\leq\xi}}\to\mod{B_\xi}.
$$
Its left adjoint   
\begin{equation}\label{EStandFunct}
\funE_\xi=\De_\xi\otimes_{B_\xi}-:\mod{B_\xi}\to \mod{H_{\leq\xi}}
\end{equation}
is called a {\em weak standardization functor}. 
By \cite[Lemma 4.11]{Kdonkin}, if $\varrho(\pi)=\xi$ then 
$$\De(\pi)\cong \funE_\xi(P_\xi(\pi))\quad \text{and} \quad \bar\De(\pi)\cong \funE_\xi(L(\pi)).$$ 
A weak standardization functor is called a {\em standardization functor} if it  is exact. This is equivalent to the requirement that $\De(\xi)$ is flat as a $B_\xi$-module.

\begin{Definition}\label{DStCat} 
{\rm 

The algebra $H$ as above is called  {\em properly stratified } (with respect to the fixed preorder $\leq$) if the following properties hold:
\begin{enumerate}
\item[{\tt (Filt)}] For every $\pi\in\Pi$, the object $K(\pi)$ has a finite $\De$-filtration with quotients of the form  $q^d \De(\si)$ for $\si>\pi$. 

\item[{\tt (Flat)}] For every  $\xi\in\Xi$, the right $B_\xi$-module $\De(\xi)$ is finitely generated and flat. 
\end{enumerate}
}
\end{Definition}

\section{Lie theoretic preliminaries}\label{SLie}
\subsection{Affine root system}\label{SSARS}
Let $\Car=(\cc_{ij})_{i,j\in I}$ be a {\em Cartan matrix} of  untwisted affine $ADE$ type, see \cite[\S 4, Table Aff 1]{Kac}. 
We have $I=\{0,1,\dots,l\},$ 
where $0$ is the affine vertex. Following \cite[\S 1.1]{Kac}, let $(\h,\Pi,\Pi^\vee)$ be a realization of the Cartan matrix $\Car$, so we have simple roots $\{\al_i\mid i\in I\}$ 
and a bilinear form $(\cdot,\cdot)$ on $\h^*$ such that 
$\cc_{ij}=2(\al_i,\al_j)/(\al_i,\al_i)$ for all $i,j\in I$. We normalize $(\cdot,\cdot)$ so that $(\al_i,\al_i)=2$ for all $i\in I$. 
Denote 
$Q_+ := \bigoplus_{i \in I} \Z_{\geq 0} \cdot \al_i$. For $\theta \in Q_+$, we write $\height(\theta)$ for the sum of its 
coefficients when expanded in terms of the $\al_i$'s.

Let $\g'=\g(\Car')$ be the finite dimensional simple Lie algebra  whose Cartan matrix $\Car'$ corresponds to the subset of vertices $I':=I\setminus\{0\}$. We denote by $W$ and $W'$ the corresponding affine and finite Weyl groups. 
Let $\Phi'$ and $\Phi$ be the root systems of $\g'$ and $\g$ respectively. Denote by $\Phi'_+$ and $\Phi_+$ the set of {\em positive}\, roots in $\Phi'$ and $\Phi$, respectively, cf. \cite[\S 1.3]{Kac}. Let 
$\de\in \Phi_+$ 
be the {\em null root}. 
We have 
$\Phi_+=\Phi_+^\im\sqcup \Phi_+^\re
$, where
$\Phi_+^\im=\{n\de\mid n\in\Z_{>0}\}$
and 
\begin{align*}
\Phi_+^\re=\{\be+n\de\mid \be\in  \Phi'_+,\ n\in\Z_{\geq 0}\}\sqcup \{-\be+n\de\mid \be\in  \Phi'_+,\ n\in\Z_{> 0}\}.
\end{align*}
Let $$p:\Phi\to\Phi'$$ be the projection which maps $n\de$ to $0$ and $\pm\be+n\de$ to $\pm\be$. For $\al\in\Phi'$ we denote by $\hat\al$ the minimal height root in $\Phi_+^\re$ with $p(\hat\al)=\al$.

We denote by $$\Psi:=\Phi_+^\re\sqcup\{\de\}$$ the set of {\em indivisible positive roots}.

\subsection{Convex preorders}\label{SSConv}
A {\em convex preorder} on $\Phi_+$ is a total preorder $\preceq$ such that for all $\be,\ga\in\Phi_+$ we have:
\begin{enumerate}
\item[$\bullet$]
If $\be\preceq \ga$ and $\be+\ga\in\Phi_+$, then $\be\preceq\be+\ga\preceq\ga$;
\item[$\bullet$]
$\be\preceq\ga$ and $\ga\preceq\be$ if and only if $\be$ and $\ga$ are imaginary.
\end{enumerate}
Thus, upon restriction to $\Psi$, the preorder $\preceq$ is a total order. 
We write $\be\prec\ga$ if $\be\preceq\ga$ but $\ga\not\preceq\be$. The set of real roots splits into two disjoint infinite sets
$$
\Phi_{\succ \de}:=\{\be\in \Phi_+\mid \be\succ\de\}\quad \text{and}\quad  
\Phi_{\prec \de}:=\{\be\in \Phi_+\mid \be\prec\de\}.
$$
A convex preorder is called {\em balanced} if $p(\Phi_{\succ \de})=\Phi'_+$. In general we have:

\begin{Lemma} \label{LpPosRoots} {\rm \cite[Lemma 3.7]{McNAff}} 
There is $w\in W'$ such that $p(\Phi_{\succ \de})=w\Phi'_+$. 
\end{Lemma}

Let $p(\Phi_{\succ \de})=w\Phi'_+$ according to the lemma. For $i\in I'$, we denote 
$$\ga_i:=w\al_i\quad\text{and}\quad \ga_i^\pm:=\widehat{\pm\ga_i}.
$$ 
Then $\ga_i^\pm\in\Phi_+^\re$, $\ga_i^++\ga_i^-=\de$, and $\ga_i^+\succ\de\succ\ga_i^-$. Note that $$\De_{\succ \de}:=\{\ga_1,\dots,\ga_l\}$$ is a base in $\Phi'$.

\begin{Lemma} \label{LSpecialOrder} {\rm \cite[Example 3.5]{McNAff}} 
Let $\De$ be any base in $\Phi'$ and $\al\in \De$. There exists a convex preorder $\preceq$ on $\Phi_+$ with the following three properties:
\begin{enumerate}
\item[{\rm (i)}] $\De_{\succ \de}=\De$;
\item[{\rm (ii)}] $\hat \al\succ\hat \al+\de\succ \hat \al+2\de\succ\dots\succ\de\dots\succ\widehat{-\al}+2\de\succ\widehat{-\al}+\de\succ\widehat{-\al}$;
\item[{\rm (iii)}] Every root in $\Phi_+^\re$, which is not of the form $\widehat{\pm\al}+n\de$, is either greater than $\hat\al$ or less than $\widehat{-\al}$.
\end{enumerate} 
\end{Lemma}

In this subsection we write $\equiv$ for $\equiv\pmod{\Z\de}$.

\begin{Lemma} \label{LDiff} 
Let $i\in I'$ and  
$\ga_i^\pm=\eta^\pm + \theta^\pm$ for some $\eta^\pm,\theta^\pm\in Q_+$. Suppose that $\eta^\pm$ is a sum of positive roots $\preceq \ga_i^\pm$, and $\theta^\pm$ is a sum of positive roots $\succeq \ga_i^\pm$. Then $\eta^-+\eta^+\neq \ga_i^-$ unless $\eta^+=\theta^-=0$.
\end{Lemma}
\begin{proof}
By assumption, $\theta^+$ is a sum of positive roots $\succeq \ga_i^+$. 
Since $\ga_i^+\succ\de$, these positive roots are in $\Phi_{\succ \de}$. So we can write $\theta^+\equiv\sum_{m=1}^l c_m\ga_m$ with coefficients $c_m\in\Z_{\geq 0}$. 
Furthermore, $\eta^-$ is a sum of positive roots less than $\ga_i^-$. As $\ga_i^-\prec \de$, these positive roots are in $\Phi_{\prec \de}$. So we can write $\eta^-\equiv-\sum_{m=1}^l d_m\ga_m$ with coefficients $d_m\in\Z_{\geq 0}$. Assume that $\eta^-+\eta^+=\ga_i^-$, in which case $\theta_-=\eta_+$. 
As $\{\de,\ga_1,\dots,\ga_l\}$ are linearly independent, we deduce that all $d_m$ and $c_m$ with $m\neq i$ are zero and $d_i+c_i=2$. 

If $d_i=1$ then $\theta^-=\eta^+\equiv 0$, which implies $\theta^-=\eta^+= 0$ by height considerations.
If $d_i=0$, we have $\eta^-\equiv 0$ so $\eta^-= 0$ by heights, which implies $\theta^-=\eta^+=\ga_i^-$ and $\theta_+=\ga_i^+ -\ga_i^-\equiv 2\ga_i$. Since $\height(\theta^+)<\height(\ga_i^+)$, we deduce that $\theta^+$ is a sum of positive roots which are strictly greater than $\ga_i^+$. As $(\theta_+,\ga_i^+)=4$, any presentation of $\theta^+$ as a sum of positive roots which are strictly greater than $\ga_i^+$ must have at least four summands. Each of these summands is a non-negative linear combination of $\ga_m$'s. This contradicts $\theta^+\equiv 2\ga_i$. 
The case $d_i=2$ is similar to the case $d_i=0$. 
\end{proof}

\begin{Lemma} \label{LMGG}
Let $n\in\Z_{>0}$ and $\de=\theta_r^-+\theta_r^+$ for $r=1,\dots,n$, with each  $\theta_r^-$ being a sum of positive roots $\preceq \de$ and each $\theta_r^+$ being a sum of positive roots $\succeq \de$.  If\, $\sum_{r=1}^n\theta_r^\pm=n\ga_i^\pm$, then $\theta_r^\pm=\ga_i^\pm$ for all $r=1,\dots n$. 
\end{Lemma}
\begin{proof}
For $1\leq r\leq n$ we have $\theta_r^\pm\equiv \pm\sum_{j=1}^lc_{r,j}^\pm\ga_j^+$ for some  $c_{r,i}^\pm\in\Z_{\geq 0}$. 
So $\pm n\ga_i\equiv\pm\sum_{r=1}^n \sum_{j=1}^lc_{r,j}^\pm\ga_j^+$. 
Now linear independence of the $\ga_j^\pm$ modulo $\Z\de$ and considerations of height imply $c_{r,j}^\pm=\de_{i,j}$ for all $r$. 
\end{proof}

\subsection{Kostant partitions and root partitions}\label{SSKostant}
Let $\theta\in Q_+$.  
A {\em Kostant partition of $\theta$} is a tuple $\xi=(x_\be)_{\be\in\Psi}$ of non-negative integers such that $\sum_{\be\in\Psi}x_\be\be=\theta$. 
If $\{\be_1\succ\dots\succ\be_r\}=\{\be\in\Psi\mid x_\be\neq 0\}$, then,  denoting $x_u:=x_{\be_u}$, we also write $\xi$ in the form 
$$\xi=(\be_1^{x_{1}},\dots,\be_r^{x_{r}}).$$ We denote by $\Xi(\theta)$ the set of all Kostant partitions of $\theta$. Denoting the left (resp. right) lexicographic order on $\Xi(\theta)$ by $\leq_l$ (resp. $\leq_r$), we will always use the {\em bilexicographic partial order} $\leq$ on $\Xi(\theta)$, i.e. 
$\xi\leq \zeta$ if and only if $\xi\leq_l \zeta$ and $\xi\geq_r \zeta$.

Let $\al\in\Psi$. By convexity, the Kostant partition $(\al)$ is  the unique minimal element in $\Xi(\al)$. 
A {\em minimal pair for $\al$} is a minimal element in $\Xi(\al)\setminus\{(\al)\}$. A minimal pair for $\al$ exists, provided $\al$ is not a simple root, which we always assume when speaking of minimal pairs for $\al$. 
Using the property (Con2) from \cite[\S3.1]{Kcusp}, it is easy to see that a minimal pair is always a Kostant partition of the form $(\be,\ga)$ for $\be,\ga\in\Psi$ with $\be>\ga$. 
A minimal pair $(\be,\ga)$ is called {\em real} if both $\be$ and $\ga$ are real.

\begin{Lemma} \label{LRealMP}  Let $\al\in\Phi_+^\re$. If $\al$ has no real minimal pair, then  $\al=\ga_i^\pm +n\de$ for some $i\in I'$ and $n\in\Z_{> 0}$, in which case  $(\ga_i^++(n-1)\de,\de)$ is a minimal pair for $\ga_i^++n\de$ and $(\de, \ga_i^-+(n-1)\de)$ is a minimal pair for $\ga_i^-+n\de$.
\end{Lemma}
\begin{proof}
The first half is \cite[Lemma 12.4]{McNAff}. For the second half, if $(\ga_i^++(n-1)\de,\de)$ is not a minimal pair for $\ga_i^++n\de$, then we would be able to write $\ga_i^++n\de=\be+\ga$ with $\be,\ga\in\Phi_{\succ \de}$. But modulo $\de$ both $\be$ and $\ga$ are positive sums of $\ga_j$'s, which leads to a contradiction. The argument for $\ga_i^-+n\de$ is similar. 
\end{proof}

If $\mu$ is a (usual integer) partition of $n$, we write $\mu\vdash n$ and $n=|\mu|$. 
An {\em $l$-multipartition} of $n$ is a tuple $\umu=(\mu^{(1)},\dots,\mu^{(l)})$ of partitions such that $|\umu|:=|\mu^{(1)}|+\dots+|\mu^{(l)}|=n$. The set of the all $l$-multipartitions of $n$ is denoted by $\Par_n$, and $\Par:=\sqcup_{n\geq 0}\Par_n$.

A {\em root partition of $\theta$} is a pair $(\xi,\umu)$, where $\xi=(x_\be)_{\be\in\Psi}\in\Xi(\theta)$ and $\umu\in\Par_{x_\de}$. 
We write $\Pi(\theta)$ for the set of root partitions of~$\theta$. 
There is a natural surjection $\rho:\Pi(\theta)\to\Xi(\theta),\ (\xi,\umu)\mapsto \xi$. The bilexicographic partial order $\leq$ on $\Xi(\theta)$ induces a partial preorder $\leq$ on $\Pi(\theta)$, i.e. $\pi\leq \si$ if and only if $\rho(\pi)\leq \rho(\si)$.

Let $(\xi,\umu)\in \Pi(\theta)$ with $\xi=(x_\be)_{\be\in\Psi}$. 
As all but finitely many integers $x_\be$ are zero, there is a finite subset
$
\be_1\succ\dots\succ\be_s\succ\de\succ\be_{-t}\succ\dots\succ \be_{-1}
$
of $\Psi$ such that $x_\be=0$ for $\be\in\Psi$ outside of this subset. Then, denoting $x_u:=x_{\be_u}$, we can 
 write any root partition of $\theta$ in the form
\begin{equation}\label{ERP}
(\xi,\umu)=(\be_1^{x_1},\dots,\be_s^{x_s},\umu,\be_{-t}^{x_{-t}},\dots,\be_{-1}^{x_{-1}}),
\end{equation}
where all $x_u\in\Z_{\geq 0}$, $\umu\in\Par$, and 
$|\umu|\de+\sum x_u\be_u=\theta.$

\section{KLR algebras}\label{SKLR}

\subsection{Presentation}\label{SSDefKLR}
In this subsection, $k$ is an arbitrary commutative unital ring.   
Define the polynomials $\{Q_{ij}(u,v)\in k[u,v]\mid i,j\in I\}$  
as follows. For the case where the Cartan matrix $\Car\neq {\tt A}_1^{(1)}$, 
choose signs $\eps_{ij}$ for all $i,j \in I$ with $\cc_{ij}
< 0$  so that $\eps_{ij}\eps_{ji} = -1$ and set 
\begin{equation*}\label{EArun}
Q_{ij}(u,v):=
\left\{
\begin{array}{ll}
0 &\hbox{if $i=j$;}\\
1 &\hbox{if $\cc_{ij}=0$;}\\
\eps_{ij}(u^{-\cc_{ij}}-v^{-\cc_{ji}}) &\hbox{if $\cc_{ij}<0$.}
\end{array}
\right.
\end{equation*}
For type $A_1^{(1)}$ we set
\begin{equation*}\label{EArun1}
Q_{ij}(u,v):=
\left\{
\begin{array}{ll}
0 &\hbox{if $i=j$;}\\
(u-v)(v-u) &\hbox{if $i\neq j$.}
\end{array}
\right.
\end{equation*}

Fix 
$\theta\in Q_+$ of height $n$. Let
$$
I^\theta=\{\bi=(i_1, \dots, i_n)\in I^n\mid \al_{i_1}+\dots+\al_{i_n}=\theta\}. 
$$ 
The symmetric group $\Si_n$ with simple transpositions $s_r:=(r,r+1)$ acts on  $I^\theta$ by place permutations. 

The {\em KLR-algebra} $R_\theta$ is an associative graded unital $k$-algebra, given by the generators
$
\{1_{\bi}\mid \bi\in I^\theta\}\cup\{y_1,\dots,y_{n}\}\cup\{\psi_1, \dots,\psi_{n-1}\}
$ 
and the following relations for all $\bi,\bj\in I^\theta$ and all admissible $r,t$:
\begin{equation*}
1_{\bi}  1_{\bj} = \de_{\bi,\bj} 1_{\bi} ,
\quad{\textstyle\sum_{\bi \in I^\theta}} 1_{\bi}  = 1;\label{R1}
\end{equation*}
\begin{equation*}\label{R2PsiY}
y_r 1_{\bi}  = 1_{\bi}  y_r;\quad y_r y_t = y_t y_r;
\end{equation*}
\begin{equation*}
\psi_r 1_{\bi}  = 1_{s_r\bi} \psi_r;\label{R2PsiE}
\end{equation*}
\begin{equation*}
(y_t\psi_r-\psi_r y_{s_r(t)})1_{\bi}  
= \de_{i_r,i_{r+1}}(\de_{t,r+1}-\de_{t,r})1_{\bi};
\label{R6}
\end{equation*}
\begin{equation*}
\psi_r^21_{\bi}  = Q_{i_r,i_{r+1}}(y_r,y_{r+1})1_{\bi} 
 \label{R4}
\end{equation*}
\begin{equation*} 
\psi_r \psi_t = \psi_t \psi_r\qquad (|r-t|>1);\label{R3Psi}
\end{equation*}
\begin{equation*}
\begin{split}
&(\psi_{r+1}\psi_{r} \psi_{r+1}-\psi_{r} \psi_{r+1} \psi_{r}) 1_{\bi}  
\\=
&
\de_{i_r,i_{r+2}}\frac{Q_{i_r,i_{r+1}}(y_{r+2},y_{r+1})-Q_{i_r,i_{r+1}}(y_r,y_{r+1})}{y_{r+2}-y_r}1_{\bi}.
\end{split}
\label{R7}
\end{equation*}
The {\em grading} on $R_\theta$ is defined by setting 
$
\deg(1_{\bi} )=0$, $\deg(y_r1_{\bi} )=2$, $\deg(\psi_r 1_{\bi} )=-(\al_{i_r},\al_{i_{r+1}}).
$

\subsection{Basic representation theory of KLR algebras}\label{SSBasicRep} 
For any $V\in\mod{R_\theta}$, its {\em formal character} $\CH V:=\sum_{\bi\in \words^\theta}(\DIM 1_{\bi} V)\cdot\bi$ is an element of $\bigoplus_{\bi\in\words^\theta}\Z((q))\cdot \bi$. We will use  the fact proved in \cite{KL1} that the formal characters of irreducible $R_\theta$-modules are linearly independent. 
Note also that $\CH(q^dV)=q^d\CH(V)$, where the first $q^d$ means the degree shift. We refer to $1_{\bi} V$ as the {\em $\bi$-weight space} of $V$ and to its vectors as {\em vectors of weight}~$\bi$.

There is an anti-automorphism $\tau:R_\theta \rightarrow R_\theta$ which fixes all standard generators. Given $V\in\mod{R_\theta}$, we denote $V^\circledast := \HOM_{k}(V, k)$ viewed as a left $R_\theta$-module via $\tau$. Note that in general $V^\circledast$ is not finitely generated as an $R_\theta$-module, but we will apply $\circledast$ only to finite dimensional modules. We have $\CH V^{\circledast}=\overline{\CH V}$, where the bar means the {\em bar-involution}, i.e. the automorphism of $\Z[q,q^{-1}]$ that swaps $q$ and $q^{-1}$ extended to $\bigoplus_{\bi\in \words^\theta}\Z[q,q^{-1}]\cdot \bi$. 

It is shown in \cite{KL1} that the algebra $R_\theta$ is Noetherian, Laurentian and Schurian. Moreover, there is always a unique choice of degree shift for every irreducible module $L$ such that $L^\circledast\cong L$.

\begin{Lemma} \label{LGeneralHead} 
Let $V\in\mod{R_\theta}$, $\bi \in I^\theta$, and $v\in 1_{\bi} V$ be a non-zero homogeneous vector with $R_\theta v=V$. Assume that there is only one irreducible $R_\theta$-module $L$ up to $\simeq$ with $1_{\bi} L\neq 0$ and $[V:L]_q\neq 0$. Then $\head V\simeq L$. 
\end{Lemma}
\begin{proof}
If $W$ is the radical of $V$ then $V/W\cong \oplus_{r}m_r(q) L_r$ for simple modules $L_r$, with $L_r\not\cong L_s$ for $r\neq s$, and multiplicities $m_r(q)\in\Z[q,q^{-1}]$. By assumptions, there exists $r$ such that $L\cong L_r$, $m_s(q)=0$ for $s\neq r$, and $v+W\in m_r(q) L_r$. Finally, $v+W$ generates $m_r(q) L_r$, so  $m_r(q)$ is of the form $q^d$. 
\end{proof}

\subsection{Induction and restriction}

Let $\theta_1,\dots,\theta_m\in Q^+$ and $\theta=\theta_1+\dots+\theta_m$. Denote $\underline{\theta}:=(\theta_1,\dots,\theta_m)$. 
Consider the set of concatenations  $$I^{\underline{\theta}}:=\{\bi^1\cdots\bi^m \mid \bi^1\in\words^{\theta_1},\dots,\bi^m\in\words^{\theta_m}\}\subseteq \words^\theta.$$ 
There is a natural (non-unital) algebra embedding 
$$R_{\underline{\theta}}=R_{\theta_1,\dots,\theta_m}:=R_{\theta_1}\otimes\dots\otimes R_{\theta_m}\to R_\theta,$$ 
which sends the unit $1_{\theta_1}\otimes\dots \otimes 1_{\theta_m}$ to the idempotent 
$ 
1_{\underline{\theta}}:=\sum_{\bi\in \words^{\underline{\theta}}}1_{\bi} \in R_\theta.
$ 
We have the 
induction functor 
$$\Ind_{\underline{\theta}}=\Ind_{\theta_1,\dots,\theta_m}:=R_\theta 1_{\underline{\theta}}\otimes_{R_{\underline{\theta}}}-:\mod{R_{\underline{\theta}}}\to\mod{R_{\theta}}.
$$ 
For $V_1\in\mod{R_{\theta_1}}, \dots, V_m\in\mod{R_{\theta_m}}$, we denote 
$$V_1\circ\dots\circ V_m:=\Ind_{\underline{\theta}} V_1\boxtimes \dots\boxtimes V_m.$$ 

We also have the restriction functors: 
$$
\Res_{\underline{\theta}} =\Res_{\theta_1,\dots,\theta_m}:= 1_{\underline{\theta}} R_{\theta}
\otimes_{R_{\underline{\theta}}} -:\mod{R_{\theta}}\rightarrow \mod{R_{\underline{\theta}}}.
$$
While the induction functor $\Ind_{\underline{\theta}}$ is left adjoint to $\Res_{\underline{\theta}}$,  its right adjoint is given by the coinduction:
$$
\Coind_{\underline{\theta}}:=\Hom_{R_{\underline{\theta}}}(1_{\underline{\theta}}R_{\theta},\,-) 
$$

From  {\rm \cite[Theorem 2.2]{LV}} and \cite[Lemma 2.21]{Kcusp} we have

\begin{Lemma} \label{LLV} 
Let $\underline{\theta}:=(\theta_1,\dots,\theta_m)\in Q_+^m$, and $V_k\in\mod{R_{\theta_k}}$ for $k=1,\dots,m$. Set $d(\underline{\theta}):=\sum_{1\leq l<k\leq m}(\theta_l,\theta_k)$. 
Then 
$$
(V_1\circ\dots\circ V_m)^\circledast\cong 
q^{d(\underline{\theta})}
V_m^\circledast\circ\dots\circ V_1^\circledast\cong 
\Coind_{\underline{\theta}}\, 
V_1^\circledast\boxtimes\dots\boxtimes V_m^\circledast.
$$
\end{Lemma}

All these functors have parabolic analogues. For example, given a family $(\theta^a_b)_{1\leq a\leq n,\ 1\leq b\leq m}$ of elements of $Q_+$, set 
$\sum_{a=1}^n\theta^{a}_b=:\theta_b$ for all $1\leq b\leq m$. Then we have obvious functors
$$
\Ind_{(\theta^{1}_1,\dots,\theta^{n}_{1});\dots;(\theta^{1}_{m},\dots,\theta^{n}_{m})}^{\theta_1;\dots;\theta_m}:
\mod{R_{\theta^{1}_1}\otimes \dots\otimes R_{\theta^{n}_{m}}}\to
\mod{R_{\theta_1}\otimes \dots\otimes R_{\theta_m}}.
$$

Given $\underline{\theta}=(\theta_1,\dots,\theta_N)\in Q_+^N$, and a permutation $x\in \Si_N$, we denote 
$
x\underline{\theta}:=(\theta_{x^{-1}(1)},\dots,\theta_{x^{-1}(N)}) 
$, 
and define the integer
$$
s(x,\underline{\theta}):=-\sum_{1\leq m<k\leq N,\ x(m)>x(k)}(\theta_m,\theta_k).
$$
There is an obvious algebra isomorphism 
$
\phi^x:R_{x\underline{\theta}}\to R_{\underline{\theta}}
$
permuting the components. Composing with this isomorphism, we get a functor
$$
\mod{R_{\underline{\theta}}}\to \mod{R_{x\underline{\theta}}},\  M\mapsto {}^{\phi^x}M.
$$
Making an additional shift, we get a functor 
$$
\mod{R_{\underline{\theta}}}\to \mod{R_{x\underline{\theta}}},\  M\mapsto {}^xM:=
q^{ s(x,\underline{\theta})}
({}^{\phi^x}M).
$$

Fix $\underline{\eta}=(\eta^1,\dots,\eta^n)\in Q_+^n$ and $\underline{\theta}=(\theta_1,\dots,\theta_m)\in Q_+^m$  
with 
$$
\eta^1+\dots+\eta^n=\theta_1+\dots+\theta_m.
$$ 
Let $\D(\underline{\theta},\underline{\eta})$ be the set of all tuples $\underline{\kappa}=(\kappa^a_b)_{1\leq a\leq n,\ 1\leq b\leq m}$ of elements of $Q_+$ such that  $\sum_{b=1}^m\kappa^{a}_b=\eta^a$ for all $1\leq a\leq n$ and $\sum_{a=1}^n\kappa^{a}_b=\theta_b$ for all $1\leq b\leq m$. 
For each $\underline{\kappa}\in \D(\underline{\theta},\underline{\eta})$, we define the permutation $x(\underline{\kappa})\in \Si_{mn}$ which maps  
$$
(\kappa^{1}_1,\dots,\kappa^{1}_{m},\kappa^{2}_1,\dots,\kappa^{2}_{m}\dots,\kappa^{n}_{1},\dots,\kappa^{n}_{m})
$$
to
$$
(\kappa^{1}_1,\dots,\kappa^{n}_{1},\kappa^{1}_2,\dots,\kappa^{n}_{2},\dots,\kappa^{1}_{m},\dots,\kappa^{n}_{m}).
$$

Let $M\in\mod{R_{\underline{\eta}}}$. We can now consider the $R_{\kappa^{1}_1,\dots,\kappa^{n}_{1};\dots;\kappa^{1}_{m},\dots,\kappa^{n}_{m}}$-module 
$${}^{x(\underline{\kappa})}\big(\Res_{\kappa^{1}_1,\dots,\kappa^{1}_{m};\dots;\kappa^{n}_{1},\dots,\kappa^{n}_{m}}^{\eta^1;\dots;\eta^n}
\,M \big).$$

We will need the following weak version of the Mackey Theorem for KLR algebras, see the proof of \cite[Proposition 2.18]{KL1}:

\begin{Theorem} \label{TMackeyKL}
Let  $M\in\mod{R_{\underline{\eta}}}$. The  $R_{\underline{\theta}}$-module $\Res_{\underline{\theta}}\,\Ind_{\underline{\eta}} M$ has a filtration with subquotients   
$$
\Ind_{\kappa^{1}_1,\dots,\kappa^{n}_{1};\dots;\kappa^{1}_{m},\dots,\kappa^{n}_{m}}^{\theta_1;\dots;\theta_m}
{}^{x(\underline{\kappa})}\big(\Res_{\kappa^{1}_1,\dots,\kappa^{1}_{m};\dots;\kappa^{n}_{1},\dots,\kappa^{n}_{m}}^{\eta^1;\dots;\eta^n}
\,M \big),
$$
one for each $\underline{\kappa}\in \D(\underline{\theta},\underline{\eta})$. 
\end{Theorem}

\subsection{Reduction modulo $p$}
The KLR algebra $R_\theta$ is defined over an arbitrary commutative unital  ring $k$, and if we need to to emphasize which $k$ we are working with, we will use the notation $R_{\theta,k}$. Likewise in the notation for modules. Let $p$ be a fixed prime number, and $F:=\Z/p\Z$ be the prime field of characteristic~$p$. We will use the $p$-modular system $(F,\O,K)$ with $F={\mathbb F}_p$, $\O=\Z_p$ and $K=\Q_p$.

Let $k=K$ or $F$, and $V_k$ be an $R_{\theta,k}$-module. An $R_{\theta,\O}$-module $V_\O$ is called an {\em $\O$-form of $V_k$} if every graded component of $V_\O$ is free of finite rank as an $\O$-module and, identifying $R_{\theta,\O}\otimes_\O k$ with $R_{\theta,k}$, we have  $V_\O\otimes_\O k\cong V_k$ as $R_{\theta,k}$-modules. 
Every $V_K\in \mod{R_{\theta,K}}$ has an $\O$-form: pick $R_{\theta,K}$-generators $v_1,\dots,v_r$ and define $V_\O:=R_{\theta,\O}\cdot v_1+\dots+R_{\theta,\O}\cdot v_1$. We always can and will pick the generators which are homogeneous  weight vectors. 
Let $V_K\in\mod{R_{\theta,K}}$ and $V_\O$ be an $\O$-form of $V_K$. 
The $R_{\theta,F}$-module $V_\O\otimes_\O F$ is called a {\em reduction modulo $p$} of $V_K$. Reduction modulo $p$ in general depends on the choice of $V_\O$. However, as explained in \cite[Lemma 4.3]{KS}, we have a generalization of the standard result for finite groups:

\begin{Lemma} \label{LRedIndep} 
If $V_{K}\in\mod{R_{\theta,K}}$ and $L_F$ is an irreducible $R_{\theta,F}$-module, then the multiplicity $[V_\O\otimes_\O F:L_F]_q$ is independent of the choice of an $\O$-form $V_\O$ of $V_K$. 
\end{Lemma}

Reduction modulo $p$ commutes with induction and restriction \cite[Lemma 4.5]{KS}:

\begin{Lemma} \label{LIndScal} 
Let $\underline{\theta}=(\theta_1,\dots,\theta_m)\in Q_+^m$, $\theta=\theta_1+\dots+\theta_m$, $V_\O\in \mod{R_{\underline{\theta};\O}}$, and $W_\O\in\mod{R_{\theta,\O}}$. Then for any $\O$-algebra $k$, there are natural isomorphisms of $R_{\theta,k}$-modules 
$$(\Ind_{\underline{\theta}} V_\O)\otimes_\O k\cong \Ind_{\underline{\theta}} (V_\O\otimes_\O k)$$ 
and of $R_{\underline{\theta};k}$-modules 
$$(\Res_{\underline{\theta}} W_\O)\otimes_\O k\cong \Res_{\underline{\theta}} (W_\O\otimes_\O k).$$ 
\end{Lemma}

In particular, reduction modulo $p$ preserves formal characters. This fact together with linear independence of  formal characters of irreducible modules has the following consequence:

\begin{Lemma} \label{LRedLinInd} 
Let $V_1,\dots, V_r$ be $R_{\theta,K}$-modules such that $\CH V_1,\dots,\CH V_r$ are linearly independent. Let $L_1,\dots, L_s$ be a complete set of composition factors of reductions modulo $p$ of the modules $V_1,\dots, V_r$. Then $s\geq r$.  
\end{Lemma}

\section{Semicuspidal modules}
\label{SSCSM}

The main goal of this section is to generalize some results of \cite{KM} and \cite{McNAff}. The paper \cite{KM} assumes that the convex order is balanced, while \cite{McNAff} assumes that $p=0$. We want to avoid both of these assumptions. 

In this section we often work with a composition $\nu=(n_1,\dots,n_h)$ of $n$, the corresponding parabolic subalgebra $$R_{\nu,\de}:=R_{n_1\de,\dots, n_h\de},$$ and the corresponding induction and restriction functors 
$$
\HCI_\nu^n:=\Ind_{n_1\de,\dots,n_h\de}\quad \text{and}\quad \HCR_\nu^n:=\Res_{n_1\de,\dots,n_h\de}.
$$

\subsection{Semicuspidal modules}\label{SSSemiCusp}
We fix a convex preorder $\preceq$ on $\Phi_+$,  an  {\em indivisible positive root}\, $\al$, and $n\in\Z_{>0}$. 
Following \cite{McNAff} (see also \cite{KRbz,McN,Kcusp,TW}) an $R_{n\al}$-module $V$ is called {\em semicuspidal}\, if $\theta,\eta\in Q_+$, $\theta+\eta=n\al$, and $\Res_{\theta,\eta}V\neq 0$ imply that $\theta$ is a sum of positive roots  $\preceq\al$ and $\eta$ is a sum of positive roots $\succeq\al$. 

Weights $\bi\in I^{n\al}$, which appear in some semicuspidal $R_{n\al}$-modules, are called {\em semicuspidal weights}. We denote by $I^{n\al}_{\noncusp}$ the set of non-semicuspidal weights. Let
$$
1_{\noncusp}
:=\sum_{\bi\in I^{n\al}_{\noncusp}}1_{\bi}.
$$
Following \cite{McNAff}, define the {\em semicuspidal algebra} 
\begin{equation}\label{ESCA}
C_{n\al}=C_{n\al,k}:=R_{n\al}/R_{n\al} 1_{\noncusp}R_{n\al}.
\end{equation}
Then the category of finitely generated semicuspidal $R_{n\al}$-modules is equivalent to the category $\mod{C_{n\al}}$.

\begin{Theorem} \label{TIrrCusp}  
Let $\al\in \Phi_+^\re$ and $n\in\Z_{>0}$. 
There is a unique up to isomorphism irreducible $\circledast$-self-dual semicuspidal $R_\al$-module. We denote it $L(\al)$. Moreover, $L(\al^n):=q^{n(n-1)/2}L(\al)^{\circ n}$ is the  unique up to isomorphism irreducible $\circledast$-self-dual semicuspidal $R_{n\al}$-module. 
\end{Theorem}
\begin{proof}
Follows from \cite[Main Theorem and Lemmas 3.3, 4.6]{Kcusp}, see also \cite{TW}.
\end{proof}

\begin{Lemma} \label{LRedSemiCuspReal}
Let $\al\in \Phi_+^\re$ and $n\in \Z_{>0}$. Then $L(\al^n)_F$ is a reduction modulo $p$ of $L(\al^n)_K$.
\end{Lemma}
\begin{proof}
See \cite[Proposition 4.9]{Kcusp} and \cite[Lemma 4.6]{KS}.
\end{proof}

In the rest of this section we work with the imaginary case trying to understand the irreducible semicuspidal $R_{n\de}$-modules.

\subsection{Minuscule imaginary modules} \label{SSMinusc}

The proof of the following lemma in \cite[Lemma 12.3]{McNAff} seems to need the assumption $p=0$, but the same result will later follow in general from Lemmas~\ref{LRedSemiCuspReal} and \ref{LLDe}(iii). Recall the definition of a minimal pair from \S\ref{SSKostant}.

\begin{Lemma} \label{LCuspMinPairRes}%
{\rm \cite[Lemma 12.3]{McNAff}} 
Assume that $p=0$. Let $\al\in\Psi$, $L\in \mod{R_\al}$ be a semicuspidal module, and $(\be,\ga)$ be a minimal pair for $\al$. Then all composition factors of $\Res_{\ga,\be}L$ are of the form $L_\ga\boxtimes L_\be$, where $L_\ga$ is an irreducible semicuspidal $R_\ga$-module and $L_\be$ is an irreducible  semicuspidal $R_\be$-module. 
\end{Lemma}

By \cite{Kcusp}, there are exactly $|\Par_1|=l$ isomorphism classes of self-dual irreducible semicuspidal $R_\de$-modules. These modules can be labeled canonically by the elements of $I'$, see  \cite{Kcusp} for balanced convex orders, and \cite{McNAff,TW} in general. We now describe the approach of \cite{McNAff}. One needs to be careful to make sure that the assumption $p=0$ made in \cite{McNAff} can be avoided. Recall the base $\De_{\succ \de}=\{\ga_1,\dots,\ga_l\}$ in $\Phi'_+$ from $\S\ref{SSConv}$ and the roots $\ga_i^\pm$. In characteristic  zero, parts (i) and (ii) of the following result are contained in \cite{McNAff}, and this will be used in the proof.

\begin{Lemma} \label{LLDe}
Let $i\in I'$. Then the module $L(\ga_i^-)\circ L(\ga_i^+)$ has a simple $\circledast$-self-dual head. Moreover, denoting this simple module by $L_{\de,i}$, we have the following:
\begin{enumerate}
\item[{\rm (i)}] The $R_\de$-module $L_{\de,i}$ is cuspdial, and $\{L_{\de,i}\mid i\in I'\}$ is a complete and irredundant system of irreducible $\circledast$-self-dual semicuspidal $R_\de$-modules. 
\item[{\rm (ii)}] $\Res_{\ga_i^-,\ga_i^+}L_{\de,i}\cong L(\ga_i^-)\boxtimes L(\ga_i^+)$, and $\Res_{\ga_j^-,\ga_j^+}L_{\de,i}=0$ if $i\neq j$. 
\item[{\rm (iii)}] Reduction modulo $p$ of $L_{\de,i,K}$ is $L_{\de,i,F}$.
\end{enumerate} 
\end{Lemma}
\begin{proof}
By Mackey's Theorem and Lemma~\ref{LDiff}, we have 
$$\Res_{\ga_i^-,\ga_i^+}(L(\ga_i^-)\circ L(\ga_i^+))\cong L(\ga_i^-)\boxtimes L(\ga_i^+).$$ 
If $L$ is a simple constituent of the head, then $L(\ga_i^-)\boxtimes L(\ga_i^+)$ appears in the socle of $\Res_{\ga_i^-,\ga_i^+} L$. Since $\Res$ is an exact functor and the multiplicity of the irreducible $\circledast$-self-dual module $L(\ga_i^-)\boxtimes L(\ga_i^+)$ in $\Res_{\ga_i^-,\ga_i^+}(L(\ga_i^-)\circ L(\ga_i^+))$ is $1$, the head is simple and self-dual. The first part of (ii) also follows.

Now, we explain that in characteristic zero, (i) and (ii) are contained in \cite{McNAff}. Indeed, (i) is \cite[Theorem 17.3]{McNAff}. To see the second part of (ii), in view of \cite[Theorem 13.1]{McNAff} and Lemma~\ref{LSpecialOrder}, we may assume that $(\ga_j^+,\ga_j^-)$ is a minimal pair for $\de$, in which case by Lemma~\ref{LCuspMinPairRes}, all composition factors of $\Res_{\ga_j^-,\ga_j^+} L_{\de,i}$ are of the form $L(\ga_j^-)\boxtimes L(\ga_j^+)$, in particular, $L(\ga_j^-)\boxtimes L(\ga_j^+)$ appears in the socle of $\Res_{\ga_j^-,\ga_j^+} L_{\de,i}$, whence $L_{\de,i}$ is a quotient of $L(\ga_j^-)\circ L(\ga_j^+)$, i.e. 
$L_{\de,i}\cong L_{\de,j}$, giving a contradiction. 

Pick $R$-forms $L(\ga_i^\pm)_R$ of $L(\ga_i^\pm)_K$. 
By Lemmas~\ref{LIndScal} and \ref{LRedSemiCuspReal}, we have that $L(\ga_i^-)_R\circ L(\ga_i^+)_R$ is an $R$-form of $L(\ga_i^-)_k\circ L(\ga_i^+)_k$ for $k=K$ or $F$. We have a surjection 
$\phi:L(\ga_i^-)_K\circ L(\ga_i^+)_K\to L_{\de,i,K}$. Let $L_{\de,i,R}:=\phi(L(\ga_i^-)_R\circ L(\ga_i^+)_R)$. Note that $L_{\de,i,R}$ is an $R$-form of $L_{\de,i,K}$. On the other hand, we have a surjection $L(\ga_i^-)_F\circ L(\ga_i^+)_F\to L_{\de,i,R}\otimes _R F$. This implies that $L_{i,\de,F}$ is a quotient of $L_{\de,i,R}\otimes _R F$. As $L_{\de,i,K}$ is semicuspidal by \cite{McNAff}, it now follows that so is $L_{\de,i,F}$. 

Let $j\neq i$. By the characteristic zero result, we have $\Res_{\ga_j^-,\ga_j^+} L_{\de,i,K}=0$. It now follows that $\Res_{\ga_j^-,\ga_j^+} L_{\de,i,F}=0$, too, which completes the proof of (ii). By (ii), we conclude that $L_{\de,i,F}\not\cong L_{\de,j,F}$. By counting simple semicuspidal $R_\de$-modules, we complete the proof of (i). 

To prove (iii), note by characters that all composition factors of $L_{\de,i,R}\otimes _R F$ are semicuspidal. Now we can conclude that $L_{\de,i,R}\otimes _R F\cong L_{\de,i,F}$ using (ii).  
\end{proof}

\begin{Remark} \label{R240616} 
{\rm 
Lemmas~\ref{LLDe}(iii) and \ref{LRedSemiCuspReal} show that the statement of Lemma~\ref{LCuspMinPairRes} holds without the assumption $p=0$. 
}
\end{Remark}

Following the terminology of \cite{Kcusp}, we call the modules $L_{\de,i}$ {\em minuscule modules}.

\subsection{Imaginary Schur-Weyl duality}
\label{SSISW}

Fix $i\in I'$. Recall the  minuscule module $L_{\de,i}$ from \S\ref{SSMinusc}. Consider the $R_{n\de}$-module 
$$M_{n,i}:=L_{\de,i}^{\circ n}$$ called the {\em imaginary tensor space of color $i$}, and the algebra 
$$\ImS_{n,i}:=R_{n\de}/\operatorname{Ann}_{R_{n\de}}(\Mde_{n,i})$$ 
called the {\em imaginary Schur algebra of color $i$}. 
We have the following result, the proof of which, given in \cite[Theorems 4.3.2, 4.4.1]{KM}, does not use the fact that the convex order is balanced. 

\begin{Theorem} \label{TSW}
Let $i\in I'$ and $n\in\Z_{>0}$. Then:
\begin{enumerate}
\item[{\rm (i)}] $\End_{R_{n\de}}(M_{n,i})^\op\cong\operatorname{end}_{R_{n\de}}(M_{n,i})^\op\cong k\Si_n$. 
\item[{\rm (ii)}] $M_{n,i}$ is a projective $\ImS_{n,i}$-module, and $M_{n,i}^\circledast\cong M_{n,i}$.
\item[{\rm (iii)}] Assume that $p>n$ or $p=0$. Then $\ImS_{n,i}$ is semisimple, $M_{n,i}$ is a projective generator over $\ImS_{n,i}$, and $\ImS_{n,i}$ is Morita equivalent to $k\Si_n$. 
\end{enumerate}
\end{Theorem}

By Theorem~\ref{TSW}, if $p=0$, the number of composition factors of $M_{n,i}$ is equal to the number of partitions of $n$. Now using reduction modulo $p$ argument involving Lemmas~\ref{LLDe}(iii) and \ref{LRedLinInd}, we deduce that the same is true in general:

\begin{Lemma} \label{LNumberCompFact} 
The number of composition factors of $M_{n,i}$, up to isomorphism and degree shift, is equal to the number of partitions of $n$. 
\end{Lemma}

Recall the roots $\ga_i^+$ and $\ga_i^-$ from \S\ref{SSConv} and the notation $\La(h,n)$ for the set compositions of $n$ with $h$ parts. 

\begin{Lemma} \label{LExtrResM}
Let $i\in I'$, $n\in\Z_{>0}$ and $(n_1,\dots,n_l)\in\La(l,n)$. We have:
\begin{enumerate}
\item[{\rm (i)}] $\Res_{n\ga_i^-,n\ga_i^+}M_{n,i}\cong L((\ga_i^-)^n)\boxtimes L((\ga_i^+)^n)$.
\item[{\rm (ii)}] $\Res_{n\ga_i^-,n\ga_i^+}(M_{n_1,1}\circ\dots\circ M_{n_l,l})=0$ unless $n_j=0$ for all $j\neq i$. 
\end{enumerate} 
\end{Lemma}
\begin{proof}
Follows using Mackey's Theorem and Lemmas~\ref{LMGG}, \ref{LLDe}(ii).
\end{proof}

For $\al\in\Phi_+^\re$, we denote by $P(\al^n)$ the projective cover of the irreducible semicuspidal module $L(\al^n)$. We will use a special projective module, which we we refer to as a {\em Gelfand-Graev module}. Note that its definition is  different from  the one in \cite{KM} even for balanced orders:
$$
\Ga_{n,i}:=P((\ga_i^-)^n)\circ P((\ga_i^+)^n).
$$

\begin{Lemma} \label{LOneDim} 
Let $i\in I'$, $n\in\Z_{>0}$ and $(n_1,\dots,n_l)\in\La(l,n)$. We have:
\begin{enumerate}
\item[{\rm (i)}] $\DIM \Hom_{R_{n\de}}(\Ga_{n,i},M_{n,i})=1$;
\item[{\rm (ii)}] $\Hom_{R_{n\de}}(\Ga_{n,i},M_{n_1,1}\circ\dots\circ M_{n_l,l})=0$ unless $n_j=0$ for all $j\neq i$. 
\end{enumerate}
\end{Lemma}
\begin{proof}
For any $R_{n\de}$-module $M$, we have 
$$\Hom_{R_{n\de}}(\Ga_{n,i},M)\cong\Hom_{R_{n\ga_i^-,n\ga_i^+}}(P((\ga_i^-)^n)\boxtimes P((\ga_i^+)^n),\Res_{n\ga_i^-,n\ga_i^+}M).
$$ 
So the result follows from Lemma~\ref{LExtrResM}. 
\end{proof}

Denote by ${\bf 1}_{\Si_n}$ the trivial (right) $k\Si_n$-module. Note that $\Hom_{R_{n\de}}(\Ga_{n,i},M_{n,i})$ is naturally a right $k\Si_n$-module, since $\Si_n$ acts on $M_{n,i}$ on the right in view of Theorem~\ref{TSW}(i). Since this module is $1$-dimensional by Lemma~\ref{LOneDim}(i), it is either the trivial or the sign module. If it happens to be the sign module, we redefine the right action of $\Si_n$ on $M_{n,i}$ by tensoring it with the sign representation. So we may assume without loss of generality that 
\begin{equation}\label{ETrivial}
\Hom_{R_{n\de}}(\Ga_{n,i},M_{n,i})\cong {\bf 1}_{\Si_n}\qquad(i\in I').
\end{equation}

Recall the notation $\HCI_\nu^n$ and $\HCR_\nu^n$ from the beginning of this section. For a composition $\nu=(n_1,\dots,n_h)\in\La(h,n)$, we define the $R_{\nu,\de}$-modules
$$
M_{\nu,i}:=M_{n_1,i}\boxtimes\dots\boxtimes M_{n_h,i},\quad \Ga_{\nu,i}:=\Ga_{n_1,i}\boxtimes\dots\boxtimes \Ga_{n_h,i},\quad\text{and}\quad \Ga^\nu_i:=\HCI_\nu^n \Ga_{\nu,i}.
$$ 
We have the parabolic analogue $\ImS_{\nu,i}$ of $\ImS_{n,i}$ defined as 
$$
\ImS_{\nu,i}:=R_{\nu,\de}/\Ann_{R_{\nu,\de}}(M_{\nu,i})\cong \ImS_{n_1,i}\otimes\dots\otimes \ImS_{n_h,i}.
$$
The functors $\HCR^n_\nu$ and $\HCI^n_\nu$ induce the functors between $\mod{\ImS_{n,i}}$ and $\mod{\ImS_{\nu,i}}$.

\begin{Lemma} 
Let $i\in I'$. Then $\HCR_\nu^n \Ga_{n,i}\cong \Ga_{\nu,i}\oplus X_i$, where $X_i$ is a projective $R_{\nu,\de}$-module with $\Hom_{R_{\nu,\de}}(X_i,M_{\nu,j})=0$ for all $j\in I'$. 
\end{Lemma}
\begin{proof}
Mackey's Theorem yields a filtration of 
$$\HCR_\nu^n \Ga_{n,i}=\Res_{n_1\de,\dots,n_h\de} \Ind_{n\ga_i^-,n\ga_i^+} P((\ga_i^-)^n)\boxtimes P((\ga_i^+)^n)$$ with projective subquotients, one of which is $\Ga_{\nu,i}$ (ignoring grading shifts for now). So we get a decomposition $\HCR_\nu^n \Ga_{n,i}\cong q^d\Ga_{\nu,i}\oplus X_i$ where $X_i$ is a projective module. It remains to notice that $$\DIM \Hom_{R_{\nu,\de}}(\HCR_{\nu,i}^n \Ga_{n,i},M_{\nu,j})=\de_{i,j},$$ which is done using adjointness of $\Res$ and $\Coind$ exactly as in the proof of \cite[5.1.3(i)]{KM}.
\end{proof}

\subsection{Divided powers}\label{SSDP}
Fix again $i\in I'$. 
Set ${\tt x}_n:=\sum_{x\in \Si_n}x$. 
Set 
$$X_{n,i}:=M_{n,i}x_n\quad \text{and}\quad 
Z_{n,i}:=\{v\in M_{n,i}\mid vx=v\ \text{for all $x\in \Si_n$}\}.
$$
Fixing a non-zero $R_{n\de}$-homomorphism $f_{n,i}:\Ga_{n,i}\to M_{n,i}$, we also set 
$$Y_{n,i}:=\im f_{n,i}\subseteq M_{n,i},$$ 
cf. Lemma~\ref{LOneDim}(i). Eventually we will prove that $Y_{n,i}=Z_{n,i}$. For now, it is only clear from (\ref{ETrivial}) that $Y_{n,i}\subseteq Z_{n,i}$.
For the proof of the following lemma see \cite[\S5.2]{KM}:

\begin{Lemma} \label{LXYZ}
We have: 
\begin{enumerate}
\item[{\rm (i)}] $X_{n,i}$ is an irreducible $R_{n\de}$-module. 
\item[{\rm (ii)}] $\soc Z_{n,i}=X_{n,i}$, and no composition factor of $Z_{n,i}/X_{n,i}$ is isomorphic to a submodule of $M_{n,i}$.
\end{enumerate}
\end{Lemma}

From now on fix $h\geq n$. For a letter $L\in\{X,Z,Y,\Gamma\}$ and a composition $\nu=(n_1,\dots,n_h)\in\La(h,n)$, we set 
$$L_{\nu,i}:=L_{n_1,i}\boxtimes\dots\boxtimes L_{n_h,i}, \quad L^\nu_i:=\HCI_\nu^n L_{\nu,i},\quad \text{and}\quad L_i(h,n):=\bigoplus_{\nu\in \La(h,n)}L^\nu_i.
$$ 
For the proof of the following lemma see \cite[\S5.3,5.5]{KM}.

\begin{Lemma} 
For $L\in\{X,Z,Y\}$, we have $\HCR^n_\nu L_{n,i}\cong L_{\nu,i}$. 
\end{Lemma}

Let $S(h,n)$ be the classical Schur algebra of \cite{Green}. For the proof of the following theorem see \cite[\S5.4]{KM}:

\begin{Theorem} \label{THowe} 
For $L\in\{X, Z\}$,  there is an algebra isomorphism 
$$\End_{R_{n\de}}(L_i(h,n))^{\op}=\operatorname{end}_{R_{n\de}}(L_i(h,n))^{\op}\cong S(h,n).$$ 
\end{Theorem}

The proof of the following lemma follows that of \cite[Lemma 5.5.3]{KM}:

\begin{Lemma} \label{LDoubleCosets}
Let $\la,\mu\in \La(h,n)$. Then 
\begin{align*}
\DIM \Hom_{R_{n\de}}(Y^\la_i,Y^\mu_i)&=\DIM \Hom_{R_{n\de}}(\Ga^\la_i,Y^\mu_i)\\
&=\DIM \Hom_{R_{n\de}}(\Ga^\la_i,Z^\mu_i)
=|\Si_\la\backslash\Si_n/\Si_\mu|.
\end{align*}
\end{Lemma}

We give a slightly simpler proof of the following result compared to  \cite[Theorem 5.5.4]{KM}:

\begin{Theorem} \label{TIHD}
Let $i\in I'$. For $h\geq n$, we have:
\begin{enumerate}
\item[{\rm (i)}] $Z_i(h,n)=\bigoplus_{\nu\in \La(h,n)}Z^\nu_i$ is a projective generator for $\ImS_{n,i}$.

\item[{\rm (ii)}] $Z_{n,i}=Y_{n,i}$.
\end{enumerate}
\end{Theorem}
\begin{proof}
(i) As $Y_{\nu,i}$ is a non-zero submodule of $Z_{\nu,i}$, it contains the simple socle $X_{\nu,i}$ of $Z_{\nu,i}$, see Lemma~\ref{LXYZ}. Applying $\HCI_\nu^n$ to the embeddings $X_{\nu,i}\subseteq  Y_{\nu,i}\subseteq  Z_{\nu,i}$, we get embeddings $X^\nu_i\subseteq Y^\nu_i\subseteq  Z^\nu_i$. By Lemma~\ref{LDoubleCosets}, 
\begin{align*}
\DIM \Hom_{R_{n\de}}(Y_i(h,n),Y_i(h,n))&=\DIM \Hom_{R_{n\de}}(\Ga_i(h,n),Y_i(h,n))
\\
&=\DIM \Hom_{R_{n\de}}(\Ga_i(h,n),Z_i(h,n))
\\
&=\sum_{\la,\mu\in \La(h,n)}|\Si_\la\backslash\Si_n/\Si_\mu|\\&=\dim S(h,n),
\end{align*}
the last equality for the dimension of the classical Schur algebra being well-known. 

In particular, this implies that $Y_i(h,n)$ is projective as an $R_{n\de}/\Ann_{R_{n\de}}(Y_i(h,n))$-module by Lemma~\ref{LSchubert}. But $M_{n,i}=Y^{(1^n)}_i$ is a summand of $Y_i(h,n)$. So $$\Ann_{R_{n\de}}(Y_i(h,n))=\Ann_{R_{n\de}}(M_{n,i}),$$ and  $Y_i(h,n)$ is a projective $\ImS_{n,i}$-module. By the classical theory \cite{Green}, the number of isomorphism classes of irreducible $S(h,n)$-modules equals  the number of partitions of $n$. By Fittings' Lemma, the number of isomorphism classes of indecomposable summands of $Y_i(h,n)$ equals  the number of isomorphism classes of irreducible modules over $\End_{R_{n\de}}(Y_i(h,n))=S(h,n)$. We now deduce from Lemma~\ref{LNumberCompFact} that $Y_i(h,n)$ is a projective generator for $\ImS_{n,i}$. 

(ii) By (i), every irreducible $\ImS_{n,i}$-modules appears in the head of the projective $R_{n\de}$-module $\Ga_i(h,n)$. 
As 
$$\DIM \Hom_{R_{n\de}}(\Ga_i(h,n),Y_i(h,n))
=\DIM \Hom_{R_{n\de}}(\Ga_i(h,n),Z_i(h,n)),
$$ 
every homomorphism from $\Ga_i(h,n)$ to $Z_i(h,n)$ has image lying in $Y_i(h,n)$, and it follows that $Y_i(h,n)=Z_i(h,n)$.  
\end{proof}

\subsection{Imaginary semicuspidal irreducible and Weyl modules}\label{SSMoritaEquivBe}
Recall that we have fixed $h\geq n$. By Theorem~\ref{TIHD}, for every $i\in I'$, we may regard $Z_i(h,n)$ as a $(\ImS_{n,i},S(h,n))$-bimodule. Then by Morita theory, we have an equivalences of categories 
\begin{align*}
\be_{h,n;i}&:\mod{S(h,n)}\to \mod{\ImS_{n,i}},\quad W\mapsto Z_i(h,n)\otimes_{S(h,n)}W.
\end{align*}

By the classical theory \cite{Green}, the Schur algebra $S(h,n)$ is quasihereditary with  irreducible module $L_{\tt cl}(\la)$ and standard modules $W_{\tt cl}(\la)$ labeled by the partitions $\la\vdash n$. Define the $\ImS_{n,i}$-modules: 
\begin{align*}
L_i(\la) &:= \be_{h,n;i}(L_{\tt cl}(\la))\\
W_i(\la) &:= \be_{h,n;i}(W_{\tt cl}(\la))
\end{align*}
so that  $\{L_i(\la)\mid\la\vdash n\}$ is a complete and irredundant family of irreducible modules over $\ImS_n=\ImS_{n,i}$ up to isomorphism and degree shift (as in \cite[Lemma 6.1.3]{KM} one checks that the definitions do not depend on $h\geq n$). By inflating, these are irreducible semicupsidal $R_{n\de}$-modules. It is easy to see that $L_i(\la)^\circledast\cong L_i(\la)$.

Now we complete a classification of the irreducible semicuspidal $R_{n\de}$-modules. To every multipartition $\umu=(\mu^{(1)},\dots,\mu^{(l)})\in\Par_n$, we associate the $R_{n\de}$-module
$$
L(\umu):=L_1(\mu^{(1)})\circ\dots\circ L_l(\mu^{(l)}).
$$

The proof of the following two results is the same as that of \cite[Theorem 5.10, Lemma 5.11]{Kcusp}.

\begin{Theorem} \label{TRedOneCol} 
Let $n\in\Z_{>0}$. Then $\{L(\umu)\mid\umu\in\Par_n\}$ is a complete and irredundant set of $\circledast$-selfdual  irreducible semicuspidal $R_{n\de}$-modules up to isomorphism.
\end{Theorem}

\begin{Proposition} \label{PRedOneCol} 
Let $n\in\Z_{>0}$, $\umu\in\Par_n$, $\nu=(n_1,\dots,n_l)\in\La(l,n)$, and $\ula\in\Par_n$ satisfy $|\la^{(i)}|=n_i$ for all $i$. Then:
\begin{enumerate}
\item[{\rm (i)}] All composition factors of $\HCR^n_\nu L(\umu)$ are of the form 
$$
(L_1(\mu^{(11)})\circ\dots\circ L_l(\mu^{(l1)}))\boxtimes\dots\boxtimes (L_1(\mu^{(1l)})\circ\dots\circ L_l(\mu^{(ll)}))
$$
with $\sum_{i=1}^l |\mu^{ij}|=n_j$ for all $j\in I'$ and $\sum_{j=1}^l |\mu^{ij}|=|\mu^{(i)}|$ for all $i\in I'$. 

\item[{\rm (ii)}]  If 
$L_1(\la^{(1)})\boxtimes\dots\boxtimes L_l(\la^{(l)})$ is a composition factor of $\HCR^n_\nu L(\umu)$, then $|\mu^{(i)}|=n_i$ for all $i$, $\ula=\umu$, and the multiplicity of this factor is $1$. 
\end{enumerate} 
\end{Proposition}

\begin{Corollary} \label{CCleverSubmodule}
Let $\nu=(n_1,\dots,n_l)\in\La(l,n)$. Then there exists a short exact sequence of $R_{\nu,\de}$-modules 
$$
0\to Z_1(h,n_1)\boxtimes \dots\boxtimes Z_l(h,n_l)\to \HCR^n_\nu (Z_1(h,n_1)\circ\dots\circ Z_l(h,n_l)) \to V\to 0,
$$
with $[V:L_1(\ula^{(1)})\boxtimes\dots\boxtimes L_l(\ula^{(l)})]=0$ unless $|\la^{(i)}|=n_i$ for all $i\in I'$. 
\end{Corollary}
\begin{proof}
By Theorem~\ref{TRedOneCol}, for $\umu\in\Par_n$, we have 
$$[Z_1(h,n_1)\circ\dots\circ Z_l(h,n_l):L(\umu)]
=\prod_{i\in I'}\de_{n_i,|\mu^{(i)}|}[Z_i(h,n_i):L(\mu^{(i)})].
$$  
Now the result follows from the second statement in Proposition~\ref{PRedOneCol}(ii). 
\end{proof}

\subsection{Proof of Theorem 1 and the first statement of Theorem 2}\label{SST1}
In this subsection we  always assume that $h\geq n$. 
The module $Z(h,n)$ from (\ref{EBigZ}) can be written as
\begin{equation}\label{EZHNNew}
Z(h,n)=\bigoplus_{(n_1,\dots,n_l)\in\La(l,n)} Z_1(h,n_1)\circ\dots\circ Z_l(h,n_l).
\end{equation}
Moreover, recall the algebra ${\bf S}(h,n)$ from (\ref{EBigSchur}):
$${\bf S}(h,n):= \bigoplus_{(n_1,\dots,n_l)\in \La(l,n)}S(h,n_1)\otimes \dots\otimes S(h,n_l).$$

By Theorem~\ref{THowe}, for every $i\in I'$, we have  $\End_{R_{n_i\de}}(Z_i(h,n_i))^\op\cong S(h,n_i)$. So, for every $(n_1,\dots,n_l)\in\La(l,n)$, we have an isomorphism of algebras
$$
S(h,n_1)\otimes\dots\otimes S(h,n_l) \iso \End_{R_{n_1\de,\dots,n_l\de}}(Z_1(h,n_1)\boxtimes\dots\boxtimes Z_l(h,n_l))^\op,
$$
which by functoriality of induction yields an embedding of algebras
$$
\iota_{n_1,\dots,n_l}:S(h,n_1)\otimes\dots\otimes S(h,n_l) \to \End_{R_{n\de}}(Z_1(h,n_1)\circ\dots\circ Z_l(h,n_l))^\op. 
$$
Using adjointness of induction and restriction and Corollary~\ref{CCleverSubmodule}, we see that this embedding is an isomorphism. Taking direct sum over all $(n_1,\dots,n_l)\in\La(l,n)$ we get an embedding of algebras
$
\iota: {\bf S}(h,n)\to \End_{R_{n\de}}(Z(h,n))^\op.
$
If $(n_1,\dots,n_l)$ and $(m_1,\dots,m_l)$ are distinct composition in $\La(l,n)$, then by Theorem~\ref{TRedOneCol}, the modules $Z_1(h,n_1)\circ\dots\circ Z_l(h,n_l)$ and $Z_1(h,m_1)\circ\dots\circ Z_l(h,m_l)$ do not have composition factors in common, and so  
$$\Hom_{R_{n\de}}(Z_1(h,n_1)\circ\dots\circ Z_l(h,n_l),
Z_1(h,m_1)\circ\dots\circ Z_l(h,m_l))=0.$$ 
It follows that:

\begin{Proposition} \label{PFirst} 
The map $\iota: {\bf S}(h,n)\to \End_{R_{n\de}}(Z(h,n))^\op$ is an isomorphism of algebras. 
\end{Proposition}

Recall from the introduction the imaginary Schur algebra  $$\ImS(h,n)=R_{n\de}/\Ann_{R_{n\de}}(Z(h,n)).$$ 
Our next goal is to prove that $Z(h,n)$ is a projective generator in $\mod{\ImS(h,n)}$. 

\begin{Lemma} \label{LGaHead} 
Let $i\in I'$ and $\ula=(\la^{(1)},\dots,\la^{(l)})\in\Par_n$. Then $$\Hom_{R_{n\de}}(\Ga_i(h,n), L(\ula))=0,$$ unless $\la^{(j)}=0$ for all $j\neq i$.
\end{Lemma}
\begin{proof}
As $\Ga_i(h,n)=\bigoplus_{\mu\in\La(h,n)}\Ga^\mu_i$, it suffices to prove $\Hom_{R_{n\de}}(\Ga^\mu_i, L(\ula))=0$ for an arbitrary $\mu=(m_1,\dots,m_h)\in\La(h,n)$. But
$$
\Hom_{R_{n\de}}(\Ga^\mu_i, L(\ula))=\Hom_{R_{n\de}}(\HCI_\mu^n\Ga_{\mu,i}, L(\ula))\cong \Hom_{R_{\mu,\de}}(\Ga_{\mu,i}, \HCR_\mu^n L(\ula)). 
$$
Since $L(\ula)$ is imaginary semicuspidal, all composition factors of $\HCR_\mu^n L(\ula)$ are of the form $L^1\boxtimes\dots\boxtimes L^h$ with $L^r$ imaginary semicuspidal $R_{m_r\de}$-module for $r=1,\dots,h$. So, by Lemma~\ref{LOneDim}(ii) and Theorem~\ref{TRedOneCol}, if $\Hom_{R_{\mu,\de}}(\Ga_{\mu,i}, \HCR_\mu^n L(\ula))\neq 0$, then $\HCR_\mu^n L(\ula)$ has a composition factor of the form $L_1(\nu^{(1)})\boxtimes\dots\boxtimes L_1(\nu^{(h)})$ with $\nu^{(r)}\vdash m_r$ for $r=1,\dots,h$. By Proposition~\ref{PRedOneCol}(i), we deduce  that $\la^{(j)}=0$ for all $j\neq i$. 
\end{proof}

Now Theorem 1 follows from the Proposition~\ref{PFirst} together with:

\begin{Proposition} \label{PZ(h,n)} 
The $\ImS(h,n)$-module $Z(h,n)$ is a projective generator in $\mod{\ImS(h,n)}$. 
\end{Proposition}
\begin{proof}
By Theorems~\ref{TIHD}(i) and \ref{TRedOneCol}, every irreducible imaginary semicuspidal $R_{n\de}$-module appears in the head of  $Z(h,n)$, and every irreducible $\ImS(h,n)$-module is imaginary semicuspidal, since the surjection $R_{n\de}\to \ImS(h,n)$ obviously factors through the semicuspidal algebra  $C_{n\de}$. So it remains to prove that $Z(h,n)$ is a projective $\ImS(h,n)$-module. 

It suffices to prove that $Z_1(h,n_1)\circ\dots\circ Z_l(h,n_l)$ is a projective $\ImS(h,n)$-module for any $\nu=(n_1,\dots,n_l)\in\La(l,n)$. We have a canonical surjection 
$$
\Ga_1(h,n_1)\circ\dots\circ \Ga_l(h,n_l)\to 
Z_1(h,n_1)\circ\dots\circ Z_l(h,n_l),
$$
and by Lemma~\ref{LSchubert}, it suffices to prove that 
\begin{align*}
&\dim \Hom_{R_{n\de}}(\Ga_1(h,n_1)\circ\dots\circ \Ga_l(h,n_l), Z_1(h,n_1)\circ\dots\circ Z_l(h,n_l))
\\
=&\dim \End_{R_{n\de}}(Z_1(h,n_1)\circ\dots\circ Z_l(h,n_l)).
\end{align*}
This is done as follows: 
\begin{align*}
&\dim \End_{R_{n\de}}(Z_1(h,n_1)\circ\dots\circ Z_l(h,n_l))
\\
=&\dim \Hom_{R_{\nu,\de}}(Z_1(h,n_1)\boxtimes\dots\boxtimes Z_l(h,n_l), \HCR^n_\nu Z_1(h,n_1)\circ\dots\circ Z_l(h,n_l))
\\
=&\dim \Hom_{R_{\nu,\de}}(Z_1(h,n_1)\boxtimes\dots\boxtimes Z_l(h,n_l),Z_1(h,n_1)\boxtimes\dots\boxtimes Z_l(h,n_l))
\\=&\prod_{i\in I'}\dim \Hom_{R_{n_i\de}}(Z_i(h,n_i),Z_i(h,n_i))
\\=&\prod_{i\in I'}\dim \End_{R_{n_i\de}}(\Ga_i(h,n_i),Z_i(h,n_i))
\\
=&\dim \Hom_{R_{\nu,\de}}(\Ga_1(h,n_1)\boxtimes\dots\boxtimes \Ga_l(h,n_l),Z_1(h,n_1)\boxtimes\dots\boxtimes Z_l(h,n_l))
\\
=&\dim \Hom_{R_{\nu,\de}}(\Ga_1(h,n_1)\boxtimes\dots\boxtimes \Ga_l(h,n_l),\HCR^n_\nu Z_1(h,n_1)\circ\dots\circ Z_l(h,n_l))
\\
=&\dim \Hom_{R_{n\de}}(\Ga_1(h,n_1)\circ\dots\circ \Ga_l(h,n_l), Z_1(h,n_1)\circ\dots\circ Z_l(h,n_l)),
\end{align*}
where we have used adjointness of induction and restriction for the first and the last equalities, Corollary~\ref{CCleverSubmodule} for the second equality, Lemma~\ref{LDoubleCosets} for the fourth equality, Corollary~\ref{CCleverSubmodule} and Lemma~\ref{LGaHead} for the sixth equality.
\end{proof}

For $\nu=(n_1,\dots,n_l)\in\La(l,n)$ we set
\begin{align*}
\ImS_{h,\nu}&:=\ImS_1(h,n_1)\otimes\dots\otimes \ImS_1(h,n_l),\\ 
Z_{h,\nu}&:=Z_1(h,n_1)\boxtimes\dots\boxtimes Z_l(h,n_l)\in\mod{\ImS_{h,\nu}}.
\end{align*}

\begin{Lemma} \label{LImIndGood} 
For $\nu\in\La(l,n)$ and $V\in\mod{\ImS_{h,\nu}}$, the $R_{n\de}$-module $\HCI_\nu^nV$ factors through $\ImS(h,n)$. 
\end{Lemma}
\begin{proof}
Since $Z_i(h,n_i)$ is a projective generator in $\mod{\ImS_i(h,n_i)}$ for every $i$, we have that  $Z_{h,\nu}$ is  a projective generator in $\mod{\ImS_{h,\nu}}$. So it suffices to check the lemma for $V=Z_{h,\nu}$, but this is clear. 
\end{proof}

Denote 
$$
\ImS_{h,n}:=\bigoplus_{\nu\in\La(l,n)}\ImS_{h,\nu},\qquad 
Z_{h,n}:=\bigoplus_{\nu\in\La(l,n)}Z_{h,\nu}.
$$
By Lemma~\ref{LImIndGood} we have a functor
$$
F_{h,n}:=\bigoplus_{\nu\in\La(l,n)}\HCI_\nu^n: \mod{\ImS_{h,n}}\to \mod{\ImS(h,n)},\ \bigoplus_{\nu\in\La(l,n)} V_\nu \mapsto \bigoplus_{\nu\in\La(l,n)} \HCI_\nu^nV_\nu 
$$
On the other hand, given $V\in\mod{\ImS(h,n)}$ and $\nu=(n_1,\dots,n_h)\in\La(h,n)$, we denote by $G_{h,\nu}$ the largest $R_{\nu,\de}$-submodule of $\HCR^n_\nu$ all of whose composition factors are of the form $L_1(\la^{(1)})\boxtimes \dots\boxtimes L_l(\la^{(l)})$ with $\la^{(i)}\vdash n_i$ for all $i\in I'$.

\begin{Lemma} \label{LG} 
For $\mu=(m_1,\dots,m_h),\nu=(n_1,\dots,n_h)\in\La(h,n)$, and $\la^{(i)}\vdash m_i$ for $i=1,\dots,l$, we have:
\begin{enumerate}
\item[{\rm (i)}] $G_{h,\nu}(L_1(\la^{(1)})\circ\dots\circ L_l(\la^{(l)}))=\de_{\nu,\mu}L_1(\la^{(1)})\boxtimes\dots\boxtimes L_l(\la^{(l)})$.
\item[{\rm (ii)}] $
G_{h,\nu}(Z_1(h,m_1)\circ\dots\circ Z_l(h,m_l))=\de_{\nu,\mu}Z_{h,\nu}.
$
\end{enumerate}
\end{Lemma}
\begin{proof}
Follows from Theorem~\ref{TRedOneCol}, Proposition~\ref{PRedOneCol}(ii), and Corollary~\ref{CCleverSubmodule}.
\end{proof}

Lemma~\ref{LG}(ii) and a projective generator argument imply that for any $V\in\mod{\ImS(h,n)}$, the module $G_{h,n}(V)$ factors through $\ImS_{h,\nu}$, so we get functors
\begin{align*}
G_{h,\nu}: \mod{\ImS(h,n)}&\to\mod{\ImS_{h,\nu}}\\
G_{h,n}: \mod{\ImS(h,n)}&\to\mod{\ImS_{h,n}},\ V\mapsto \bigoplus_{\nu\in\La(h,n)}G_{h,\nu}(V).
\end{align*}

\begin{Theorem} \label{TColorEquiv} 
The functors $F_{h,n}$ and $G_{h,n}$ are quasi-inverse equivalences between $\mod{\ImS_{h,n}}$ and $\mod{\ImS(h,n)}$. 
\end{Theorem}
\begin{proof}
The exact functor $F_{h,n}$ applied to the projective generator $Z_{h,n}$ yields the projective generator $Z(h,n)$ and vice versa. Moreover, the functors establish a bijection between the simples. The result follows. 
\end{proof}

By Propositions~\ref{PZ(h,n)} and \ref{PFirst}, we have a Morita equivalence functor
$$
\al(h,n)=\Hom_{\ImS(h,n)}(Z(h,n),-): \mod{\ImS(h,n)}\to\mod{{\bf S}(h,n)}.
$$
On the other hand, we also have Morita equivalence functors
\begin{align*}
\al_{h,\nu}&:=\Hom_{\ImS_{h,\nu}}(Z_{h,\nu},-): \mod{\ImS_{h,\nu}}\to\mod{(S(h,n_1)\otimes\dots\otimes S(h,n_l))},
\\
\al_{h,n}&:=\Hom_{\ImS_{h,n}}(Z_{h,n},-): \mod{\ImS_{h,n}}\to\mod{{\bf S}(h,n)}.
\end{align*}
We have the following diagram of categories and functors: 
$$
\begin{diagram}
\node{}\node{\mod{\ImS_{h,n}}}\arrow{sw,t}{F_{h,n}}\arrow{s,r}{\al_{h,n}} \\
\node{\mod{\ImS(h,n)}} \arrow{e,t}{\al(h,n)}\node{\mod{{\bf S}(h,n)}}
\end{diagram}
$$
Note that $\al_{h,n}\cong\bigoplus_{\nu\in\La(l,n)}\al_{h,\nu}$, so 
the first statement in Theorem 2 now follows from Theorem~\ref{TRedOneCol} and the following 

\begin{Proposition} \label{PComm}
We have a natural isomorphism of functors $ \al(h,n)\circ F_{h,n}\cong \al_{h,n}$. 
\end{Proposition}
\begin{proof}
This follows from the natural isomorphisms
$$
\Hom_{\ImS(h,n)}(Z(h,n),F_{h,n} V)\cong \Hom_{\ImS(h,n)}(F_{h,n}(Z_{h,n}),F_{h,n} V)\cong \Hom_{\ImS_{h,n}}(Z_{h,n},V)
$$
for any $V\in\mod{\ImS_{h,n}}$. 
\end{proof}

\section{Stratifying KLR algebras}\label{SStratKLR}

Throughout the section $\al\in\Psi$, $\theta\in Q_+$ and $\pi\in\Pi(\theta)$.

\subsection{Semicuspidal standard modules}

For real $\al$, we denote by $\De(\al^n)$ the projective cover of $L(\al^n)$ in the category $\mod{C_{n\al}}$. 
We also denote by $\De(\umu)$ the projective cover of $L(\umu)$ in the category $\mod{C_{n\de}}$. Sometimes, we will also use a special notation  $\De_{\de,i}$ for the projective cover of $L_{\de,i}$ in $\mod{C_\de}$, in other words $\De_{\de,i}=\De(\umu(i))$, where $\umu(i)$ is the multipartition of $1$ with the only non-empty component $\umu(i)^{(i)}=(1)$.

\begin{Lemma} \label{LExt1} 
Let $\al\in\Psi$ and $V\in\mod{C_{n\al}}$. Denote $\De:=\De(\al^n)$ if $\al$ is real, and $\De:=\De(\umu)$ for any $\umu\in\Par_n$ if $\al=\de$. Then 
$\Ext^1_{R_{n\al}}(\De,V)=0$. 
\end{Lemma}
\begin{proof}
Any extension of $\De$ by $V$ belongs to $\mod{C_{n\al}}$. Since $\De$ is a projective object in $\mod{C_{n\al}}$, the extension has to split. 
\end{proof}

\begin{Lemma} \label{LStandCuspidal}
Let $\al\in\Phi_+^\re$, and $n=n_1+\dots+n_a$ for $n_1,\dots,n_a\in\Z_{\geq 0}$. Then:
\begin{enumerate}
\item[{\rm (i)}] $\De(\al)^{\circ n}\cong q^{n(n-1)/2}[n]!\,\De(\al^n)$. 
\item[{\rm (ii)}] $\Res_{n_1\al,\dots,n_a\al}\De(\al^n)\cong \De(\al^{n_1})\boxtimes\dots\boxtimes \De(\al^{n_a})$. 
\end{enumerate} 
\end{Lemma}
\begin{proof}
(i) All composition factors of $\De(\al)^{\circ n}$ are of the form $L(\al^n)$, so it is an $C_{n\al}$-module. We claim that this $C_{n\al}$-module is projective. It suffices to prove that $\Ext^1_{C_{n\al}}(\De(\al)^{\circ n},L(\al^n))=0$, which follows from $\Ext^1_{R_{n\al}}(\De(\al)^{\circ n},L(\al^n))$. But 
$$\Ext^1_{R_{n\al}}(\De(\al)^{\circ n},L(\al^n))\cong \Ext^1_{R_{\al,\dots,\al}}(\De(\al)^{\boxtimes n},\Res_{\al,\dots,\al}L(\al^n)).$$ 
Now, 
\begin{equation}\label{EResIndRealCusp}
\Res_{\al,\dots,\al}L(\al^n)\cong [n]! L(\al)^{\boxtimes n},
\end{equation} 
cf. \cite[Lemma 2.11]{BKM}, so the claim follows from the K\"unneth formula and Lemma~\ref{LExt1}. 

It follows from the previous paragraph that $\De(\al)^{\circ n}\cong m(q)\,\De(\al^n)$ for some $m(q)\in \Z[q,q^{1}]$. To prove that  $m(q)=[n]!$ it suffices to observe using (\ref{EResIndRealCusp}) that $\DIM\Hom_{R_{n\al}}(\De(\al)^{\circ n},L(\al^n))=[n]!$. 

(ii) follows from (i) and the computation of $\Res_{n_1\al,\dots,n_a\al}(\De(\al)^{\circ n})$, which is performed using Mackey's Theorem and convexity. 
\end{proof}

\subsection{Standard modules}\label{SSStandardModules}
To a Kostant partition $\xi=(\be_1^{x_{1}},\dots,\be_r^{x_{r}})\in\Xi(\theta)$ we associate a parabolic subalgebra 
$$R_\xi:=R_{x_1\be_1}\otimes\dots\otimes R_{x_r\be_r}\subseteq R_\theta$$
and the corresponding functors 
\begin{equation}\label{EResRhoPi}
\Res_\xi: \mod{R_\theta}\to \mod{R_\xi}\quad\text{and}\quad \Ind_\xi,\ \Coind_\xi:\mod{R_\xi}\to \mod{R_\theta}. 
\end{equation}

For every $\pi=(\xi,\umu)\in\Pi(\theta)$ as in (\ref{ERP}), we define  the {\em proper standard module} 
\begin{equation}\label{EPrStand}
\bar\De(\pi)=L(\be_1^{x_1}) \circ \dots\circ L(\be_s^{x_s})\circ L(\umu)\circ L(\be_{-t}^{x_{-t}})\circ \dots\circ L(\be_{-1}^{x_{-1}})=\Ind_{\xi}L_\pi,
\end{equation}
and the {\em standard module} 
\begin{equation}\label{EStand}
\De(\pi)=\De(\be_1^{x_1}) \circ \dots\circ \De(\be_s^{x_s})\circ \De(\umu)\circ \De(\be_{-t}^{x_{-t}})\circ \dots\circ \De(\be_{-1}^{x_{-1}})=\Ind_{\xi}\De_\pi,
\end{equation}
where we have used the notation 
\begin{align*}
L_\pi&:=L(\be_1^{x_1}) \boxtimes \dots\boxtimes L(\be_s^{x_s})\boxtimes L(\umu)\boxtimes L(\be_{-t}^{x_{-t}})\boxtimes \dots\boxtimes L(\be_{-1}^{x_{-1}}),
\\
\De_\pi&:=\De(\be_1^{x_1}) \boxtimes \dots\boxtimes \De(\be_s^{x_s})\boxtimes \De(\umu)\boxtimes \De(\be_{-t}^{x_{-t}})\boxtimes \dots\boxtimes \De(\be_{-1}^{x_{-1}})
\end{align*}
for modules over the parabolic subalgebra
$
R_{\xi}.
$
In Lemma~\ref{LDAgree} we will show that these definitions agree with general definitions from \S\ref{SSStObStr}. Define also 
\begin{equation}\label{ENablaDef}
\bar\nabla(\pi):=\Coind_{\xi} L_\pi\cong \bar\De(\pi)^\circledast \qquad(\pi\in\Pi(\theta)),
\end{equation}
the isomorphism coming from Lemma~\ref{LLV}.

\begin{Theorem} \label{THeadIrr} {\rm \cite{Kcusp}} 
Let $\theta\in Q_+$. We have: 
\begin{enumerate}
\item[{\rm (i)}] For every $\pi\in\Pi(\theta)$, the module  
$
\bar\De(\pi)
$ has simple head; denote it $L(\pi)$. 

\item[{\rm (ii)}] $\{L(\pi)\mid \pi\in \Pi(\theta)\}$ is a complete and irredundant system of irreducible $R_\theta$-modules up to isomorphism and degree shift.

\item[{\rm (iii)}] For every $\pi\in\Pi(\theta)$, we have $L(\pi)^\circledast\cong L(\pi)$.  

\item[{\rm (iv)}] Then in the Grothendieck group $[\mod{R_{\theta}}]$, we have $[\bar\De(\pi)]=[L(\pi)]+\sum_{\si<\pi}d_{\pi,\si}[L(\si)]$ for some $d_{\pi,\si}\in\Z[q,q^{-1}]$ (which depend on $p$). 

\item[{\rm (v)}] For all $\pi,\si\in\Pi(\theta)$, we have that $\Res_{\rho(\pi)}L(\pi)\cong L_{\pi}$ and $\Res_{\rho(\si)}L(\pi)\neq 0$ implies $\si\leq \pi$.  

\end{enumerate}
\end{Theorem}

\begin{Corollary} \label{CRestr} 
Let $\theta\in Q_+$ and $\pi,\si\in \Pi(\theta)$. 
\begin{enumerate}
\item[{\rm (i)}] $\Res_{\rho(\si)}\bar\De(\pi)\neq 0$ implies $\si\leq \pi$, and   $\Res_{\rho(\pi)}\bar\De(\pi)\cong L_\pi$. 
\item[{\rm (ii)}] $\Res_{\rho(\si)}\bar\nabla(\pi)\neq 0$ implies $\si\leq \pi$ and $\Res_{\rho(\pi)}\bar\nabla(\pi)\cong L_\pi$.
\item[{\rm (iii)}] $\Res_{\rho(\si)}\De(\pi)\neq 0$ implies $\si\leq \pi$, and   $\Res_{\rho(\pi)}\De(\pi)\cong \De_\pi$. 
\end{enumerate}
\end{Corollary}
\begin{proof}
If $\Res_{\rho(\si)}\bar\De(\pi)\neq 0$, then $\Res_{\rho(\si)}L(\pi')\neq 0$ for some composition factor $L(\pi')$ of $\bar\De(\pi)$. So, using Theorem~\ref{THeadIrr}(v), we get $\si\leq\pi'\leq \pi$. The rest of (i) follows from the exactness of $\Res$ and Theorem~\ref{THeadIrr}(iv),(v). The proofs of (ii) and (iii) are similar. 
\end{proof}

\begin{Proposition} \label{PDeNa}
Let $\theta\in Q_+$, $\pi,\si\in \Pi(\theta)$, and $m\in\Z_{\geq 0}$. Then
$$
\Ext^m_{R_\theta}(\De(\pi),\bar\nabla(\si))=0
$$
unless $\rho(\pi)=\rho(\si)$. Moreover, if $\rho(\pi)=\rho(\si)$, then  $\Ext^1_{R_\theta}(\De(\pi),\bar\nabla(\si))=0$ and $\DIM\Hom_{R_\theta}(\De(\pi),\bar\nabla(\si))=\de_{\si,\pi}$. In particular, $\head \De(\pi)\cong L(\pi)$.
\end{Proposition}
\begin{proof}
The proof follows that of \cite[Proposition 24.3]{McNAff}. 
By adjointness of $\Coind$ and $\Res$, we have  
$$
\Ext^m_{R_\theta}(\De(\pi),\bar\nabla(\si))\cong \Ext^m_{R_{\rho(\si)}}(\Res_{\rho(\si)}\De(\pi),L_\si)
$$
By Corollary~\ref{CRestr}(iii), $\Res_{\rho(\si)}\De(\pi)\neq 0$ implies $\si\leq \pi$. On the other hand, by adjointness of $\Ind$ and $\Res$, we have 
$$
\Ext^m_{R_\theta}(\De(\pi),\bar\nabla(\si))\cong \Ext^m_{R_{\rho(\pi)}}(\De_\pi,\Res_{\rho(\pi)}\bar\nabla(\si)).
$$
By Corollary~\ref{CRestr}(ii), $\Res_{\rho(\pi)}\bar\nabla(\si)\neq 0$ implies $\pi\leq \si$. So we are reduced to  $\rho(\pi)=\rho(\si)$, in which case, using Corollary~\ref{CRestr}(iii), we have
$$
\Ext^m_{R_\theta}(\De(\pi),\bar\nabla(\si))\cong \Ext^m_{R_{\rho(\si)}}(\Res_{\rho(\si)}\De(\pi),L_\si)\cong \Ext^m_{R_{\rho(\si)}}(\De_\pi,L_\si).
$$
Now, the result follows from K\"unneth formula and Lemma~\ref{LExt1}.
\end{proof}

\begin{Lemma} \label{LDAgree} 
Let $\theta\in Q_+$ and $\pi\in\Pi(\theta)$. 
\begin{enumerate}
\item[{\rm (i)}] $\De(\pi)$ is the largest quotient of $P(\pi)$ all of whose composition factors $L(\si)$ satisfy $\si\leq \pi$.
\item[{\rm (ii)}] $\bar\De(\pi)$ is the largest quotient of $P(\pi)$ which has $L(\pi)$ with multiplicity $1$ and such that all its other composition factors $L(\si)$ satisfy $\si< \pi$. 
\item[{\rm (iii)}] Let $I(\pi)$ denote the injective hull of $L(\pi)$ in the category of all graded $R_\theta$-modules. Then $\bar\nabla (\pi)$ is the largest submodule of $I(\pi)$ which has $L(\pi)$ with multiplicity $1$ and all its other composition factors $L(\si)$ satisfy $\si< \pi$. 
\end{enumerate}
\end{Lemma}
\begin{proof}
(i) Since $\head \De(\pi)\cong L(\pi)$, we have a short exact sequence $$0\to X\to P(\pi)\to \De(\pi)\to 0,$$ and it suffices to prove that $\Hom_{R_\theta}(X,L(\si))=0$ if $\si\leq \pi$. Using the long exact sequence which arises by applying $\Hom_{R_\theta}(-,L(\si))$ to the short exact sequence, we have to prove $\Ext^1_{R_\theta}(\De(\pi),L(\si))=0$ for $\si\leq \pi$. But  
$$
\Ext^1_{R_\theta}(\De(\pi),L(\si))\cong\Ext^1_{R_{\rho(\pi)}}(\De_\pi,\Res_{\rho(\pi)}L(\si)).
$$
In view of Theorem~\ref{THeadIrr}(v), we may assume that $\rho(\pi)=\rho(\si)$, in which case $\Res_{\rho(\pi)}L(\si)\cong L_\si$. Now, the result follows from K\"unneth formula and Lemma~\ref{LExt1}. 

(ii) In view of Theorem~\ref{THeadIrr}(i), we have a short exact sequence $0\to X\to P(\pi)\to \bar\De(\pi)\to 0$, and it suffices to prove that $\Hom_{R_\theta}(X,L(\si))=0$ if $\si< \pi$. Using the long exact sequence which arises by applying $\Hom_{R_\theta}(-,L(\si))$ to the short exact sequence, we have to prove $\Ext^1_{R_\theta}(\bar\De(\pi),L(\si))=0$ for $\si< \pi$. But  
$$
\Ext^1_{R_\theta}(\bar\De(\pi),L(\si))\cong\Ext^1_{R_{\rho(\pi)}}(L_\pi,\Res_{\rho(\pi)}L(\si)).
$$
An application of Theorem~\ref{THeadIrr}(v) completes the proof. 

(iii) In this proof only, we will work in the larger category of all graded $R_\theta$-modules. By (\ref{ENablaDef}), $\soc \bar\nabla(\pi)\cong L(\pi)$, so there is a short exact sequence $0\to \bar\nabla(\pi)\to I(\pi)\to X\to 0$, and it suffices to prove that $X$ does not have a submodule, all of whose irreducible subfactors are $\simeq L(\si)$ with $\si<\pi$. So it suffices to prove that $X$ does not have a {\em finitely generated} submodule $Y$, all of whose composition factors are $\simeq L(\si)$ with $\si<\pi$. 
Otherwise apply $\Hom_{R_\theta}(Y,-)$ to the short exact sequence
to get an exact sequence 
$$
\Hom_{R_\theta}(Y,I(\pi))\to\Hom_{R_\theta}(Y,X)\to \Ext^1_{R_\theta}(Y,\bar\nabla(\pi))\to 0.
$$ 
Note that the middle term of this sequence is non-zero, while the first term is zero since the socle of $I(\pi)$ is $L(\pi)$. Finally, the third term is zero. Indeed,
$$
\Ext^1_{R_\theta}(Y,\bar\nabla(\pi))=\Ext^1_{R_\theta}(Y,\Coind_{\rho(\pi)}L_\pi)\cong \Ext^1_{R_{\rho(\pi)}}(\Res_{\rho(\pi)}Y,L_\pi)=0,
$$
since $\Res_{\rho(\pi)}Y=0$ in view of Theorem~\ref{THeadIrr}(v). This  a contradiction. 
\end{proof}

\subsection{Standardization functor}\label{SSStandardFunctor}
We now want to check the condition {\tt (Flat)} from Definition~\ref{DStCat}, which guarantees existence of standardization functor.

\begin{Proposition}
Let $\pi,\si\in \Pi(\theta)$ satisfy $\rho(\pi)=\rho(\si)=:\xi$. 
 Then the natural map
$
\Hom_{R_\xi}(\De_\pi,\De_\si)\to \Hom_{R_\theta}(\De(\pi),\De(\si))$
is an isomorphism. 
\end{Proposition}
\begin{proof}
By adjointness, we have 
$\Hom_{R_\theta}(\De(\pi),\De(\si))
\cong \Hom_{R_\xi}(\De_\pi,\Res_\xi \De(\si)).
$
By Corollary~\ref{CRestr}(iii), 
$
\Res_\xi\De(\si)\cong \De_\si,
$
and the result follows.
\end{proof}

\begin{Corollary} \label{CEndSplits} 
Let $\xi\in\Xi(\theta)$, $\De(\xi)=\bigoplus_{\pi\in\rho^{-1}(\xi)}\De(\pi)$, and 
$
\De_\xi:=\bigoplus_{\pi\in\rho^{-1}(\xi)}\De_\pi.
$
Then the natural map 
$
\End_{R_\xi}(\De_\xi)\to \End_{R_\theta}(\De(\xi))$
is an isomorphism of algebras. 
\end{Corollary}

\begin{Theorem} \label{TFlatness} 
Let $\theta\in Q_+$, $\xi\in\Xi(\theta)$, $\De(\xi)=\bigoplus_{\pi\in\rho^{-1}(\xi)}\De(\pi)$, and $B_\xi:=\End_{R_\theta}(\De(\xi))^{\op}$. Then, as a right $B_\xi$-module, $\De(\xi)$ is finitely generated projective, in particular, finitely generated flat.  
\end{Theorem}
\begin{proof}
We write $\xi$ in the form 
$
\xi=(\be_1^{x_1},\dots,\be_r^{x_r})$ 
for 
$\be_1\succ\dots \succ \be_r$. Note that 
$
\End_{R_\xi}(\De_\xi)^\op\cong B_{\be_1^{x_1}}\otimes \dots\otimes B_{\be_r^{x_r}}.
$ 
So by Corollary~\ref{CEndSplits}, we have 
$
B_\xi\cong B_{\be_1^{x_1}}\otimes \dots\otimes B_{\be_r^{x_r}}.
$
Moreover, each $\De(\be_m^{x_m})$ is a projective generator in the category $\mod{C_{x_m\be_m}}$. So, by Morita theory, $\De(\be_m^{x_m})$ is finitely generated projective as a right module over its  endomorphism algebra $B_{\be_m^{x_m}}$. It follows that $\De_\xi$ is finitely generated projective as a right module over its  endomorphism algebra $B_\xi$. Finally since $R_\theta$ is free of finite rank over $R_\xi$, it follows that $\De(\xi)=\Ind_\xi\De_\xi$ is finitely generated projective over $B_\xi$. 
\end{proof}

We have established the property {\tt (Flat)} from Definition~\ref{DStCat}. The property {\tt (Filt)} 
is more difficult to check. We are missing the equality $\Ext^2_{R_\theta}(\De(\pi),\bar\nabla(\si))=0$ if $\rho(\pi)=\rho(\si)$, which is needed for standard arguments as in \cite[Theorem 3.13]{BKM} yielding a $\De$-filtration on $P(\pi)$. So we will have to proceed in a round about way using reduction modulo $p$ and the results of McNamara \cite{McNAff} who has established the result in characteristic zero. 
For now, note using Proposition~\ref{PDeNa} and the K\"unneth formula, that it suffices to prove the following for all $n\in\Z_{>0}$: 
$$
\Ext^2_{R_{n\al}}(\De(\al^n),L(\al^n))=\Ext^2_{R_{n\de}}(\De(\ula),L(\umu))=0\qquad (\al\in\Phi_+^\re,\ \ula,\umu\in\Par_n).
$$

\subsection{Boundedness} 
Let $\theta = \sum_{i\in I}a_i\alpha_i$ and $n=\height(\theta)$.  Recalling that $I=\{0,1,\dots,\l\}$, pick a permutation $(i_0, \dots, i_l)$ of $(0, \dots, l)$ with $a_{i_0}>0$, and define $\bi := i_0^{a_{i_0}} \cdots i_l^{a_{i_l}} \in I^\theta$.  Then the stabilizer of $\bi$ in $S_n$ is the standard parabolic subgroup 
$S_{\bi} := S_{a_{i_0}} \times \cdots \times S _{a_{i_l}}.$  
Let $S^{\bi}$ be a set of left coset representatives for $S_n/S_{\bi}$.  Then by \cite[Theorem 2.9]{KL1} or \cite[Proposition 3.9]{Ro}, the element
\begin{equation} \label{central element}
z=z_{\bi} := \sum_{w \in S^{\bi}} (y_{w(1)} + \cdots +y_{w(a_{i_1})})1_{w \cdot \bi}
\end{equation}
is central of degree $2$ in $R_\theta$.  Let $R_{\theta}'$ be the subalgebra of $R_\theta$ generated by
\begin{equation*}
\{1_{\bi} \mid  \bi \in I^\theta \} \cup \{\psi_r \mid  1 \leq r < n \} \cup \{ y_r-y_{r+1} \mid  1 \leq r < n \}.
\end{equation*}
The restrictions from $R_\theta$ to $R_\theta'$ of modules $L(\pi),L_{\de,i},\De(\pi)$, etc. are denotes $L'(\pi),L'_{\de,i},\De'(\pi)$, etc.

\begin{Lemma} \label{LH'}  
We have 
\begin{enumerate}
\item[{\rm (i)}] $\{(y_1-y_2)^{m_1} \cdots (y_{n-1}-y_{n})^{m_{n-1}} \tau_w 1_{\bi} \mid  m_r \in \mathbb{Z}_{\geq 0}, w \in S_n, \bi\in I^\theta\}
$ 
is a basis for $R_{\theta}'$.
\item[{\rm (ii)}] If $a_{i_0}\cdot 1_k\neq 0$ in $k$, then 
there is an algebra isomorphism
$
R_\theta \cong R_\theta' \otimes k[z].
$
\end{enumerate}
\end{Lemma}
\begin{proof}
The proof is given in \cite[Lemma 3.1]{BKOld}, see also  \cite[Lemma 3.10]{KS}.
\end{proof}

For $\theta\in\Phi_+\setminus\{n\cdot\de\mid p| n\}$, and in particular for $\theta\in\Psi$, there always exists an  index $i_0$ with $a_{i_0} \cdot 1_k \neq 0$. We always make this choice. Then:

\begin{Corollary} 
For $\al\in\Psi$, we have $
R_\al \cong R_\al' \otimes k[z].
$
\end{Corollary}

Let $\al\in\Psi$, and $L$ be an irreducible $R_\al$-module. Then $z$ acts as zero on $L$, so the restriction $L'$ is an irreducible $R_\al'$-module by the corollary. For $\al\in\Phi_+$, we consider the module over $R_\al=R_\al'\otimes k[z]$: 
\begin{equation}\label{EDeltaTilde}
\tilde\De(\al):=L'(\al)\otimes k[z].
\end{equation}
Eventually we will prove that $\tilde\De(\al)\cong \De(\al)$.

As we have pointed out in Remark~\ref{R240616}, the statement of Lemma~\ref{LCuspMinPairRes} holds without the assumption $p=0$. This statement and Lemma~\ref{LLDe}(ii) is all that is needed for the argument of \cite[Theorem 15.5]{McNAff} to go through, so we get:

\begin{Lemma} \label{LBounded} {\rm \cite[Theorem 15.5]{McNAff}} 
Let $\al\in\Psi$. Then dimension of the graded components $\dim(C_\al)_d$ are bounded as a function of $d$. 
\end{Lemma}

Note that $\tilde \De(\al)\in\mod{C_\al}$ and $F[z]$ acts on $\tilde\De(\al)$ freely, so the restriction of the natural surjection $\phi:R_\al\to C_\al$ to $F[z]$ is injective, and its image gives us a central subalgebra $F[z]\subseteq C_\al$. 
Every projective $C_\al$-module is free over the subalgebra $F[z]$, and by Lemma~\ref{LBounded}, it has to be free of finite rank. Moreover, we can write $C_\al=C_\al'\otimes F[z]$ for the finite dimensional algebra $C_\al':=\phi(R_\al')$.
The same argument works for $C_\de$. Thus:

\begin{Corollary} 
Let $\al\in\Psi$. Every  standard $R_\al$-module is free of finite rank upon restriction to the subalgebra $F[z]$. Moreover, we can represent $C_\al$ as a tensor product of algebra $C_\al\cong C_\al'\otimes F[z]$ with finite dimensional $C_\al'$. 
\end{Corollary}

It is now clear that $\De(\al)\cong P'(\al)\otimes F[z]$ and $\De_{\de,i}\cong P'_{\de,i}\otimes F[z]$, where $P'(\al)$ is the projective cover of $L'(\al)$ in $\mod{C_\al'}$ and $P'_{\de,i}$ is the projective cover of $L'_{\de,i}$ in $\mod{C_\de'}$. 
The following result in characteristic zero is obtained in \cite{McNAff}:

\begin{Lemma} \label{LExactSequencesDelta}
Let $\al\in\Psi$ and $i\in I'$.
\begin{enumerate}
\item[{\rm (i)}] If $\al\in\Phi_+$ and $(\be,\ga)$ is a real minimal pair for $\al$, then there exists a short exact sequence
$$
0\to q\,\De(\be)\circ \De(\ga)\to \De(\ga)\circ\De(\be)\to \De(\al)\to 0.
$$
\item[{\rm (ii)}] If $n>1$ and $\al=\ga_i^\pm+n\de$, then, setting $\be^\pm:=\ga_i^\pm+(n-1)\de$, there exist  short exact sequences of the form
\begin{align*}
0\to \De(\be^+)\circ \De_{\de,i}\to \De_{\de,i}\circ\De(\be^+)\to & (q+q^{-1})\De(\ga_i^++n\de)\to 0,
\\
0\to \De_{\de,i}\circ \De(\be^-)\to \De(\be^-)\circ\De_{\de,i}\to & (q+q^{-1})\De(\ga_i^--n\de)\to 0.
\end{align*}

\item[{\rm (iii)}] If $\al=\de$, then there exists a short exact sequence
$$
0\to q^2\,\De(\ga_i^+)\circ \De(\ga_i^-)\to \De(\ga_i^-)\circ\De(\ga_i^+)\to \De_{\de,i}\to 0.
$$
\end{enumerate}
\end{Lemma}
\begin{proof}
(i) Lemma~\ref{LBounded} and the central subalgebra $F[z]\subseteq C_\al$ are the main ingredients in the proof of \cite[Lemma 16.1]{McNAff}, which now goes through to give the short exact sequence 
$
0\to q\,\De(\be)\circ \De(\ga)\to \De(\ga)\circ\De(\be)\to Q\to 0,
$
where $Q$ is a projective $C_{\al}$-module. 
To prove that $Q\cong \De(\al)$ it suffices to prove that $\DIM \Hom_{R_\al}(Q,L(\al))=1$. Applying $\Hom_{R_\al}(-,L(\al))$ to the short exact sequence and observing that $\Hom_{R_\al}(\De(\be)\circ \De(\ga),L(\al))=0$ by semicuspidality of $L(\al)$, we see that it suffices to prove that $\DIM \Hom_{R_\al}(\De(\ga)\circ\De(\be),L(\al))=1$. By adjointness, this dimension equals the multiplicity $[\Res_{\ga,\be}:L(\ga)\boxtimes L(\be)]_q$. In view of Lemma~\ref{LRedSemiCuspReal}, this multiplicity is independent of the characteristic of the ground field. Since the result is true in characteristic zero by \cite{McNAff}, we deduce that it also holds in positive characteristic. 

(ii) is proved analogously to (i). 

(iii) In view of \cite[Theorem 13.1]{McNAff} and Lemma~\ref{LSpecialOrder}, we may assume that $(\ga_j^+,\ga_j^-)$ is a minimal pair for $\de$. As in (i), we have a short exact sequence 
$$
0\to q^2\,\De(\ga_i^+)\circ \De(\ga_i^-)\to \De(\ga_i^-)\circ\De(\ga_i^+)\to Q\to 0,
$$
where $Q$ is a projective $C_{\de}$-module. 
To prove that $Q\cong \De_{\de,i}$ it suffices to prove that $\DIM \Hom_{R_\de}(Q,L_{\de,j})=\de_{i,j}$, which follows by applying $\Hom_{R_\al}(-,L_{\de,j})$ to the short exact sequence and observing that $$\Hom_{R_\de}(\De(\ga_i^+)\circ \De(\ga_i^-),L_{\de,j})=0$$ by semicuspidality of $L_{\de,j}$, while 
$$\DIM \Hom_{R_\de}(\De(\ga_i^-)\circ \De(\ga_i^+),L_{\de,j})=\de_{i,j}$$ by Lemma~\ref{LLDe}(ii). 
\end{proof}

\subsection{Stratification}
Recall from the end of \S\ref{SSStandardFunctor} that to prove that $R_\al$ is 
properly stratified we need some $\Ext$-result. In this subsection we prove the missing result under an explicit restriction on $p$. Again, we follow \cite{McNAff} closely.

\begin{Lemma} \label{Lp>nProjGen}
Let $\De_\de:=\bigoplus_{i\in I'}\De_{\de,i}$. Then $\De_\de^{\circ n}$ is a projective $C_{n\de}$-module. Moreover, if $p>n$ or $p=0$, then $\De_\de^{\circ n}$ is a projective generator in $\mod{C_{n\de}}$. 
\end{Lemma}
\begin{proof}
To prove that $\De_\de^{\circ n}$ is projective in $\mod{C_{n\de}}$, it suffices to show that $\Ext^1_{C_{n\de}}(\De_\de^{\circ n},L)=0$ for any irreducible $C_{n\de}$-module $L$, which would follow from $
\Ext^1_{R_{n\de}}(\De_\de^{\circ n},L)=0$. 
But the latter $\Ext$-group is isomorphic to 
$$\Ext^1_{R_{\de,\dots,\de}}(\De_\de\boxtimes\dots\boxtimes \De_\de,\Res_{\de,\dots,\de}L),
$$ 
which is indeed trivial by K\"unneth formula, since all composition factors of $\Res_{\de,\dots,\de}L$ are of the form $L_1\boxtimes \dots\boxtimes L_n$ with each $L_r$ semicuspidal. 

To show that $\De_\de^{\circ n}$ is a projective generator, it now suffices to show that $\DIM\Hom_{R_{n\de}}(\De_\de^{\circ n},L(\umu))\neq 0$ for any $\umu\in\Par_n$, which from Theorems~\ref{TSW}(iii) and \ref{TRedOneCol}. 
\end{proof}

\begin{Theorem} \label{TExt} 
Let $\al\in\Psi$ and $n\in\Z_{>0}$.
\begin{enumerate}
\item[{\rm (i)}] Let $\al=\de$. If $p>n$ or $p=0$, then for all $\ula,\umu\in\Par_n$, we have $\Ext^m_{R_{n\de}}(\De(\ula),L(\umu))=0$ for all $m>0$.
\item[{\rm (ii)}] If $\al$ is real, then $\Ext^m_{R_{n\al}}(\De(\al^n),L(\al^n))=0$ for all $m>0$.  
\end{enumerate}
\end{Theorem}
\begin{proof}
(i) By Lemma~\ref{Lp>nProjGen}, $\De_{\de}^{\circ n}$ is a projective generator in $\mod{C_{n\de}}$, so it suffices to prove that $\Ext^m_{R_{n\de}}(\De_{\de}^{\circ n},L(\umu))=0$ for all $\umu\in\Par_n$. The last $\Ext$ group is isomorphic to 
$$\Ext^m_{R_{\de,\dots,\de}}(\De_{\de}^{\boxtimes n},\Res_{\de,\dots,\de}L(\umu)).$$ 
All composition factors of $\Res_{\de,\dots,\de}L(\umu)$ are of the form $L_1\boxtimes \dots\boxtimes L_n$ with each $L_r\in\mod{C_\de}$, so by the K\"unneth formula, we may assume that $n=1$, i.e. we need to prove $\Ext^m_{R_\de}(\De_{\de,i},L_{\de,j})=0$ for all $i,j\in I'$ and $m>0$. But this follows by applying $\Hom_{R_\de}(-,L_{\de,j})$ to the short exact sequence in Lemma~\ref{LExactSequencesDelta}(iii), using  Lemma~\ref{LLDe}(ii) and induction on the height. 

(ii) In view of Theorem~\ref{TIrrCusp}(i) and Lemma~\ref{LStandCuspidal}(i), we may assume that $n=1$. 
To prove $\Ext^m_{R_\al}(\De(\al),L(\al))=0$, we apply $\Hom_{R_\al}(-,L(\al))$ to the short exact sequence in Lemma~\ref{LExactSequencesDelta}(i),(ii), and use (i) and induction on height.
\end{proof}

Taking into account the results of \S\S\ref{SSStandardModules},\ref{SSStandardFunctor}, we now have:

\begin{Corollary} 
Let $\theta=\sum_{i\in I}n_i\al_i\in Q_+$ and assume that $p>\min\{n_i\mid i\in I\}$. For any convex preorder on $\Phi_+$, the algebra $R_\al$ is properly stratified with standard modules $\{\De(\pi)\mid \pi\in\Pi(\theta)\}$ and proper standard modules $\{\bar\De(\pi)\mid \pi\in\Pi(\theta)\}$. 
\end{Corollary}

\section{Reduction modulo $p$ of irreducible and standard modules}\label{SRed}

\subsection{Reduction modulo $p$ of irreducible modules}

We already know from Lemma~\ref{LRedSemiCuspReal} that reduction modulo $p$ of a real semicuspidal module $L(\al^n)_K$ is $L(\al^n)_F$. We now look at reductions modulo $p$ of some  imaginary semicuspidal modules. For  $\la,\mu\vdash n$, we denote by 
$$d^p_{\tt cl}(\la,\mu):=[W_{\tt cl}(\la):L_{\tt cl}(\la)]$$ 
the decomposition numbers for the classical Schur algebra $S(n,n)$ in characteristic $p$. It is known that $d^p_{\tt cl}(\la,\la)=1$ and $d^p_{\tt cl}(\la,\mu)=0$ unless $\mu\unlhd \la$ in the dominance order. 
For $\ula,\umu\in\Par_n$, we define
$$d^p(\ula,\umu):=
\left\{
\begin{array}{ll}
\prod_{i\in I'}d^p_{\tt cl}(\la^{(i)},\mu^{(i)}) &\hbox{if $|\la^{(i)}|=|\mu^{(i)}|$ for all $i\in I'$,}\\
0 &\hbox{otherwise.}
\end{array}
\right.
$$
Again, $d^p(\ula,\ula)=1$ and $d^p(\ula,\umu)=0$ unless $\umu\unlhd \ula$, which means by definition that $\mu^{(i)}\unlhd \la^{(i)}$ for all $i\in I'$. 

\begin{Lemma} \label{LImRedModP} 
Let $i\in I'$ and $\la,\mu\vdash n$. Then $W_i(\la)_F$ is reduction modulo of $W_i(\la)_K=L_i(\la)_K$. In particular, $[L_i(\la)_\O\otimes F: L_i(\mu)_F]_q=d^p_{\tt cl}(\la,\mu)$. 
\end{Lemma}
\begin{proof}
The first statement is proved exactly as \cite[Theorem 6.4.3]{KM}. The second statement now follows by the Morita equivalence $\be_n$ from \S\ref{SSMoritaEquivBe}. 
\end{proof}

\begin{Lemma} \label{LDecNumbersGeneralizedSchur} 
Let $\ula,\umu\in\Par_n$. Then $L(\ula)_\O\otimes_\O F$ is  semicuspidal, and $[L(\ula)_\O\otimes_\O F:L(\umu)_F]_q=d^p(\ula,\umu)$. 
\end{Lemma}
\begin{proof}
Induction and reduction modulo $p$ commute by Lemma~\ref{LIndScal}, so the result follows from Lemma~\ref{LImRedModP} and Theorem~\ref{TRedOneCol}.
\end{proof}

\begin{Corollary} \label{CRedLImp>n} 
For $\umu\in\Par_n$ and $p>n$, reduction modulo $p$ of $L(\umu)_K$ is $L(\umu)_F$. 
\end{Corollary}

Let $\theta\in Q_+$ and $\pi\in\Pi(\theta)$. 
Let $1_F\in R_{\theta,F}$ be a primitive idempotent such that $R_{\theta,F}1_F\cong P(\pi)_F$. By an argument in \cite[\S4.1]{KS}, there is an idempotent $1_\O\in R_{\theta,\O}$ with $1_F=1_\O\otimes 1$. Let $P(\pi)_\O:=R_{\theta,\O}1_\O$. Extending scalars to $K$ we get a projective $R_{\theta,\O}$-module $P(\pi)_\O\otimes_\O K$. So we can decompose it as a direct sum of some projective indecomposable modules $P(\si)_K$.

\begin{Lemma} \label{LAdjNumbers} 
Let $\ula\in\Par_n$ and $\pi=(\xi,\ula)\in\Pi(\theta)$. Then in the Grothendieck group $[\mod{R_{\theta,F}}]$, we have 
$$
[L(\pi)_\O\otimes_\O F]=[L(\pi)_F]+\sum_{\umu\lhd \ula}d^p(\ula,\umu)[L((\xi,\umu))_F]+\sum_{\si<\pi}a_{\pi,\si}[L(\si)_F]
$$
for some bar-invariant Laurent polynomials $a_{\pi,\si}\in\Z[q,q^{-1}]$. Moreover, 
$$
P(\pi)_\O\otimes_\O K\cong P(\pi)_K\oplus \bigoplus _{\umu\rhd\ula}d^p(\umu,\ula)P((\xi,\umu))_K\oplus\bigoplus_{\si>\pi}a_{\si,\pi}P(\si)_K.
$$
\end{Lemma}
\begin{proof}
Similar to the proof of \cite[Lemma 4.8]{KS}, but using Lemma~\ref{LDecNumbersGeneralizedSchur}.
\end{proof}

\begin{Corollary} \label{CCFReductionDelta} 
All composition factors $L(\si)_F$ of a reduction modulo $p$ of $\De(\pi)_K$ satisfy $\si\leq \pi$.
\end{Corollary}

\subsection{Reduction modulo $p$ of standard modules}
The proof of following result uses an idea from \cite{Thompson}.

\begin{Lemma} \label{LThompson} 
Let $\pi\in \Pi(\theta)$. Then $\De(\pi)_F$ contains a submodule $M$ such that $\De(\pi)_F/M$ is a reduction modulo $p$ of $\De(\pi)_K$. 
\end{Lemma}
\begin{proof}
By Lemma~\ref{LAdjNumbers}, we can decompose
$P(\pi)_\O\otimes_\O K\cong P(\pi)_K\oplus Q$ for some $R_{\theta,K}$-module $Q$. Since $\De(\pi)_K$ is a quotient of $P(\pi)_K$, there is an $R_{\theta,K}$-submodule $V_K\subseteq P(\pi)_\O\otimes_\O K$ with $P(\pi)_\O\otimes_\O K/V_K\cong \De(\pi)_K$. Let $V_\O=V_K\cap P(\pi)_\O$, where we consider $P(\pi)_\O$ as an $\O$-submodule of $P(\pi)_\O\otimes_\O K$ in a natural way. Note that $V_\O$ is a pure $R_{\theta,\O}$-invariant sublattice in $P(\pi)_\O$ and $P(\pi)_\O/V_\O$ is an $\O$-form of $\De(\pi)_K$. So 
$(P(\pi)_\O/V_\O)\otimes_\O F$, which is a reduction modulo $p$ of $\De(\pi)_K$,
is a quotient 
of $P(\pi)_F$. By Corollary~\ref{CCFReductionDelta}, all composition factors $L(\si)_F$ of $(P(\pi)_\O/V_\O)\otimes_\O F$ satisfy $\si\leq \pi$, so by definition of $\De(\pi)_F$ as the largest quotient of $P(\pi)_F$ with such composition factors,  $(P(\pi)_\O/V_\O)\otimes_\O F$ is a quotient of $\De(\pi)_F$.
\end{proof}

Let $\al\in\Phi_+^\re$ and $n\in\Z_{>0}$. 
We have a semicuspidal standard module $\De(\al^n)_K$.  Pick a generator $v\in\De(\al^n)_K$ which is a homogeneous weight vector. 
Consider the $R_{n\al,\O}$-invariant lattice 
$\De(\al^n)_\O:=R_{n\al,\O}\cdot v,$ and the  reduction $\De(\al^n)_\O\otimes_\O F$. 

By Lemma~\ref{LLDe}, we have 
$\Res_{\ga_i^-,\ga_i^+}L_{\de,i}\cong L(\ga_i^-)\boxtimes L(\ga_i^+)$ 
and $\Res_{\ga_i^-,\ga_i^+}L_{\de,j}=0$ for $j\neq i$. So, picking a weight $\bj^\pm$ of $L(\ga_i^\pm)$, we have a weight $\bj^i:=\bj^-\bj^+$ of $L_{\de,i}$ such that $1_{\bj^i}L_{\de,j}=0$ for all $j\neq i$. Pick a homogeneous generator $v\in\De_{\de,i,K}$ {\em of weight $\bj^i$}. 
Consider the invariant lattice 
$\De_{\de,i,\O}:=R_{\de,\O}\cdot v$ and the  reduction $\De_{\de,i,\O}\otimes_\O F$.

\begin{Lemma} \label{LRedDelta}
We have 
\begin{enumerate}
\item[{\rm (i)}] $\De(\al^n)_\O\otimes_\O F$ is a semicuspidal $R_{n\al,F}$-module with simple head $L(\al^n)_F$, and so it is a quotient of $\De(\al^n)_F$. 
\item[{\rm (ii)}] $\De_{\de,i,\O}\otimes_\O F$ is a semicuspidal $R_{\de,F}$-module with simple head $L_{\de,i,F}$, and so it is a quotient of $\De_{\de,i,F}$. 
\end{enumerate}
\end{Lemma}
\begin{proof}
By Lemma~\ref{LRedSemiCuspReal}, $L(\al^n)_\O\otimes_\O F\cong L(\al^n)_F$ is irreducible, so all composition factors of $\De(\al^n)_\O\otimes_\O F$ are isomorphic to $L(\al^n)_F$, i.e. this module is semicuspidal. By Lemma~\ref{LLDe}(iii), we see similarly that $\De_{\de,i,\O}\otimes_\O F$ is also semicuspidal. 
In both situations, $v\otimes 1\in \De_\O\otimes_\O F$ is a cyclic generator of $\De_\O\otimes_\O F$, and it remains to apply Lemma~\ref{LGeneralHead}. 
\end{proof}

Now we can prove a stronger result:

\begin{Theorem} \label{TReductionsofCuspidal|De} 
Let $\al\in\Phi_+^\re$ and $i\in I'$. Then $\De(\al)_F\cong \De(\al)_{\O}\otimes_\O F$ and $\De_{\de,i,F}\cong \De_{\de,i,\O}\otimes_\O F$.
\end{Theorem}
\begin{proof}
Apply induction on $\height(\al)$. The base being clear, and the inductive step is obtained from Lemmas~\ref{LRedDelta} and \ref{LExactSequencesDelta} by character considerations. 
\end{proof}

\begin{Corollary} 
If $\al\in\Phi_+^\re$, then $\De(\al)\cong \tilde\Delta(\al)$ and $\End_{R_{\al}}(\De(\al))\cong F[z]$. 
\end{Corollary}
\begin{proof}
Since $L'(\al)$ is irreducible, we deduce by adjointness that $\tilde\De(\al)$ has simple head, whence it is a quotient of $\De(\al)$. Now compare the characters using \cite[Theorem 18.3]{McNAff} in characteristic zero and Theorem~\ref{TReductionsofCuspidal|De}.
\end{proof}

\begin{Corollary} \label{CSymPol} 
If $\al\in\Phi_+^\re$ and $n\in\Z_{>0}$, then $\De(\al^n)_F\cong \De(\al^n)_\O\otimes_\O F$ and $\End_{R_{\al}}(\De(\al^n))\cong F[z_1,\dots,z_n]^{\Si_n}$.
\end{Corollary}
\begin{proof}
The first statement follows from Lemmas~\ref{LRedDelta}(i),   \ref{LStandCuspidal} and Theorem~\ref{TReductionsofCuspidal|De} by induction on $n$. The second statement then follows using the fact that it is true in characteristic zero \cite{McNAff}. 
\end{proof}

We can now prove that certain cuspidal algebras $C_\al$ are `defined over integers'. 

\begin{Corollary} 
Let $\al\in\Phi_+$ and $n\in\Z_{>0}$. Then $C_{n\al,\O}$ and $C_{\de,\O}$ are free over $\O$, with  $C_{n\al,k}\cong C_{n\al,\O}\otimes_\O k$ and $C_{\de,k}\cong C_{\de,\O}\otimes_\O k$ for  and $k=F$ or $K$.
\end{Corollary}
\begin{proof}
We explain the argument for $C_\de$, the argument for $C_{n\al}$ being similar. The isomorphisms $C_{\de,k}\cong C_{\de,\O}\otimes_\O k$ are clear, and it suffices to prove that $\DIM C_{\de,K}=\DIM C_{\de,F}$. But $\DIM C_{\de,k}=\sum_{i\in I'}(\DIM L_{\de,i})(\DIM \De_{\de,i})$, which, as we have now proved, is the same for $k=K$ and $F$. 
\end{proof}

We conjecture that a similar statement is true in general. The part which remains open is:

\begin{Conjecture}
Let $n\in\Z_{>0}$ and $k=F$ or $K$. Then $C_{n\de,\O}$ is free over $\O$ and $C_{n\de,k}\cong C_{n\de,\O}\otimes_\O k$. 
\end{Conjecture}

The only difficult thing here is to show that $C_{n\de,\O}$ has no $p$-torsion. The following result implies that $C_{n\de,\O}$ at least has no  $p$-torsion if $p>n$.

\begin{Lemma} 
Let $n\in\Z_{>0}$, $\umu\in\Par_n$, and $p>n$. Then $\De(\umu)_F$ is a reduction modulo $p$ of $\De(\umu)_K$. 
\end{Lemma}
\begin{proof}
Working over $k=F$ or $K$, by Lemma~\ref{Lp>nProjGen}, $\De_{\de}^{\circ n}$ is a projective generator in $\mod{C_{n\de}}$. 
So, we can decompose $\De_{\de}^{\circ n}=\bigoplus_{\umu\in\Par_n}m(\umu)\De(\umu)_k$ with {\em non-zero} multiplicities $m(\umu)$, which a priori might depend on $k$. Moreover, 
\begin{align*}
m(\umu)&=\DIM \Hom_{R_{n\de}}(\De_{\de}^{\circ n},L(\umu))=\DIM \Hom_{R_{\de,\dots,\de}}(\De_{\de}^{\boxtimes n},\Res_{\de,\dots,\de}L(\umu))
\\
&=\DIM \Hom_{R_{\de,\dots,\de}}((\bigoplus_{i\in I'}\De_{\de,i})^{\boxtimes n},\Res_{\de,\dots,\de}L(\umu))
\\
&=\sum_{i_1,\dots,i_n\in I'}\DIM \Hom_{R_{\de,\dots,\de}}(\De_{\de,i_1}\boxtimes \dots\boxtimes \De_{\de,i_n},\Res_{\de,\dots,\de}L(\umu))
\\
&=\sum_{i_1,\dots,i_n\in I'}[\Res_{\de,\dots,\de}L(\umu)):
L_{\de,i_1}\boxtimes \dots\boxtimes L_{\de,i_n}]_q.
\end{align*}
The last expression is independent of $k$ by Corollary~\ref{CRedLImp>n}. It now follows from Lemma~\ref{LThompson} by a character argument that $\CH \De(\umu)_F=\CH\De(\umu)_K$ and  that $\De(\umu)_F$ is a reduction modulo $p$ of $\De(\umu)_K$. 
\end{proof}

\end{document}